\documentclass[11pt]{article}
\usepackage[a4paper,left=2cm,right=2cm,top=2.5cm,bottom=2.5cm]{geometry}
\usepackage{mathtools,amsmath,amsfonts,amssymb,amsthm}
\usepackage[ruled,vlined,linesnumbered]{algorithm2e}
\usepackage{booktabs,longtable,array}
\usepackage{enumitem}
\usepackage{tikz}
\usetikzlibrary{arrows.meta, positioning}
\usepackage[colorlinks=true,linkcolor=blue,citecolor=green,urlcolor=blue,pagebackref=true]{hyperref}
\def\N{\mathbb{N}}

\newtheorem{theorem}{Theorem}[section]
\newtheorem{definition}[theorem]{Definition}
\newtheorem{proposition}[theorem]{Proposition}
\newtheorem{corollary}[theorem]{Corollary}
\newtheorem{lemma}[theorem]{Lemma}
\newtheorem{remark}[theorem]{Remark}
\newtheorem{example}[theorem]{Example}
\numberwithin{equation}{section}

\SetNlSty{textnormal}{}{}
\SetFuncSty{textsc}
\SetKwInput{KwIn}{Input}
\SetKwInput{KwOut}{Output}
\SetKwFunction{FData}{Data}
\SetKwProg{Fn}{function}{}{}
\SetKwBlock{Loop}{loop}{}
\SetKwComment{Comment}{$\triangleright$\ }{}
\SetCommentSty{textnormal}

\newenvironment{algo}[2]{%
  \par\addvspace{\medskipamount}%
  \setcounter{algocf}{\value{theorem}}%
  \begin{algorithm}[H]\caption{#1}\label{#2}%
}{%
  \end{algorithm}%
  \setcounter{theorem}{\value{algocf}}%
  \par\addvspace{\medskipamount}}

\newcommand{\CP}{\mathcal C(\mathcal P)}
\newcommand{\rk}{\mathrm{rank}}
\newcommand{\Ap}{\mathrm{Ap}}
\newcommand{\PF}{\mathrm{PF}}
\newcommand{\SG}{\mathrm{SG}}
\newcommand{\msg}{\mathrm{msg}}
\newcommand{\PFs}{\mathrm{PF}^{*}}
\newcommand{\NSPF}{\texttt{Numerical\-Semigroups\-With\-Pseudo\-Frobenius\-Numbers}}

\title{The covariety of numerical semigroups with some prescribed pseudo-Frobenius numbers}
\author{W.G. Hern\'andez-Yanes and J.C. Rosales}
\date{}

\begin{document}

\maketitle

\begin{abstract}
If $S$ is a numerical semigroup, then we denote by $m(S)$, $F(S)$ and $\mathrm{PF}(S)$ its multiplicity, its Frobenius number and its set of pseudo-Frobenius numbers, respectively. Let $\mathcal P$ be a finite nonempty set of positive integers. In this paper we prove that $\mathcal C(\mathcal P)=\{S\mid S \text{ is a numerical semigroup}, F(S)=\max(\mathcal P) \text{ and } \mathcal P\subseteq\mathrm{PF}(S)\}$ is a covariety. The associated tree provides algorithms to compute $\mathcal C(\mathcal P)$, the elements of $\mathcal C(\mathcal P)$ with a prescribed type, the elements of $\mathcal C(\mathcal P)$ with a prescribed genus, and the set $\{S\mid S \text{ is a numerical semigroup and } \mathrm{PF}(S)=\mathcal P\}$. We show that every element of $\mathcal C(\mathcal P)$ has a unique minimal system of generators relative to $\mathcal C(\mathcal P)$, and we call its cardinality the $\mathcal C(\mathcal P)$-rank. Finally, we characterize the elements of $\mathcal C(\mathcal P)$ with $\mathcal C(\mathcal P)$-rank one.

\medskip


\smallskip

\end{abstract}

\section{Introduction}\label{sec:intro}

Let $\mathbb Z$ be the set of integers and $\mathbb N=\{x\in\mathbb Z \mid x\geq 0 \}$. A \textit{submonoid} of $(\mathbb N,+)$ is a subset of $\mathbb N$ closed under addition and containing 0. A \textit{numerical semigroup} is a submonoid $S$ of $(\mathbb N,+)$ such that $\mathbb N\setminus S$ is finite. If $S$ is a numerical semigroup, then $m(S)=\min(S\setminus\{0\})$, $F(S)=\max(\mathbb Z\setminus S)$ and $g(S)=\mathrm{Card}(\mathbb N\setminus S)$, where $\mathrm{Card}(X)$ denotes the cardinality of a set $X$, are three important invariants of $S$, which we call the \textit{multiplicity}, \textit{Frobenius number} and \textit{genus} of $S$, respectively.

If $\mathcal X$ is a nonempty subset of $\mathbb N$, then we denote by $\langle \mathcal X\rangle$ the submonoid of $(\mathbb N,+)$ generated by $\mathcal X$, that is, $\langle\mathcal X\rangle=\{\lambda_1x_1+\cdots+\lambda_nx_n \mid n\in\mathbb N\setminus\{0\}, \{x_1,\ldots,x_n\}\subseteq\mathcal X,\{\lambda_1,\ldots,\lambda_n\}\subseteq\mathbb N\}$. In \cite[Lemma 2.1]{libro} it is shown that $\langle \mathcal X\rangle$ is a numerical semigroup if and only if $\gcd(\mathcal X)=1$. If $M$ is a submonoid of $(\mathbb N,+)$ and $M=\langle\mathcal X\rangle$, then we say that $\mathcal X$ is a \textit{system of generators of} $M$. If, in addition, $M\ne\langle\mathcal Y\rangle$ for all $\mathcal Y\subsetneq\mathcal X$, then we say that $\mathcal X$ is a \textit{minimal system of generators of} $M$. In \cite[Corollary 2.8]{libro} it is shown that every submonoid of $(\mathbb N,+)$ has a unique minimal system of generators and that such a system is finite. We denote by $\msg(M)$ the minimal system of generators of $M$. The cardinality of $\msg(M)$ is the \textit{embedding dimension} of $M$ and we denote it by $e(M)$. If $S$ is a numerical semigroup, then $e(S)\leq m(S)$ \cite[Proposition 2.10]{libro}; when $e(S)=m(S)$ we say that $S$ has \textit{maximal embedding dimension}, or that $S$ is a \textit{MED-semigroup} (see \cite[Chapter 2]{libro} and \cite{med}).

Following the terminology introduced in \cite{JPAA}, an integer $z$ is a \textit{pseudo-Frobenius number} of a numerical semigroup $S$ if $z\notin S$ and $z+s\in S$ for all $s\in S\setminus\{0\}$. We denote by $\mathrm{PF}(S)$ the set formed by all pseudo-Frobenius numbers of $S$. The cardinality of $\mathrm{PF}(S)$ is an important invariant of $S$ (see \cite{froberg} and \cite{barucci}) called the \textit{type} of $S$, and we denote it by $t(S)$. It is clear that $\mathrm{PF}(\mathbb N)=\{-1\}$ and that $\mathrm{PF}(S)\subseteq\mathbb N\setminus\{0\}$ if $S$ is a numerical semigroup such that $S\ne\mathbb N$. In \cite[Corollary 2.23]{libro} it is shown that if $S$ is a numerical semigroup and $S\ne\mathbb N$, then $t(S)\leq m(S)-1$. It is clear that $F(S)\in\mathrm{PF}(S)$.

If $S$ is a numerical semigroup and $n\in S\setminus\{0\}$, then the \textit{Ap\'ery set} of $n$ in $S$, named after \cite{apery}, is
\[
\Ap(S,n)=\{s\in S\mid s-n\notin S\}.
\]
It has exactly $n$ elements, one in each congruence class modulo $n$ \cite[Lemma 2.4]{libro}, and every element of $S$ is the sum of an element of $\Ap(S,n)$ and a multiple of $n$.

If $\mathcal P$ is a finite nonempty subset of $\mathbb{N}\setminus\{0\}$, then we denote by $$\mathcal C(\mathcal P)=\{S \mid S\text{ is a numerical semigroup}, \mathcal P\subseteq\mathrm{PF}(S),F(S)=\max(\mathcal P)\}.$$ The elements of $\mathcal C(\mathcal P)$ will be called \textit{$\mathcal C(\mathcal P)$-semigroups}. Our aim in this work is to study the set $\mathcal C(\mathcal P)$. In particular, we will be interested in providing an algorithm that allows us to compute all the elements of $\mathcal C(\mathcal P)$. Following the terminology introduced in \cite{covariedades}, a \textit{covariety} is a nonempty family $\mathcal C$ of numerical semigroups satisfying the following conditions:

\begin{enumerate}
    \item There exists $\Delta(\mathcal C)=\min(\mathcal C)$ (with respect to set inclusion).
    \item If $\{S,T\}\subseteq\mathcal C$, then $S\cap T\in\mathcal C$.
    \item If $S\in\mathcal C$ and $S\ne\Delta(\mathcal C)$, then $S\setminus\{m(S)\}\in\mathcal C$.
\end{enumerate}

In Section~\ref{sec:tree}, we prove that $\mathcal C(\mathcal P)$ is a covariety and $\Delta(\mathcal C(\mathcal P))=\{0,\max(\mathcal P)+1,\rightarrow\}$. The results of \cite{covariedades} allow us to organize the elements of $\mathcal C(\mathcal P)$ in the form of a tree. We characterize the children of an arbitrary vertex of this tree, and this allows us, in Section~\ref{sec:algorithm}, to provide an algorithm that enables us to compute all the elements of $\mathcal C(\mathcal P)$. In particular, we obtain an algorithm (alternative to the one given in \cite{delgado}) to compute $\{S \mid \text{$S$ is a numerical semigroup and } \mathrm{PF}(S)=\mathcal P\}$. Moreover, this algorithm allows us to address more combinatorial questions in order to find numerical semigroups with prescribed properties. For instance, it enables the construction of numerical semigroups with a prescribed set $\mathcal P$ of pseudo-Frobenius numbers and a prescribed type.

The tree has more structure than the one inherited from a general covariety. The level of a semigroup $S$ in the tree is $\max(\mathcal P)-g(S)$, and the set of pseudo-Frobenius numbers decreases strictly along every edge. Combining this with inequalities of Fr\"oberg, Gottlieb and H\"aggkvist \cite{froberg} and of Nari \cite{nari}, we show in Section~\ref{sec:applications} that the semigroups of type $k$ in $\CP$ lie in a narrow band of levels, which bounds the height of the tree. We also find, for each vertex $S$, a set of pseudo-Frobenius numbers of $S$ that remain pseudo-Frobenius numbers in every descendant of $S$; this allows the algorithms to discard whole subtrees. In particular, the numerical semigroups $S$ with $\PF(S)=\mathcal P$ are obtained without building the whole tree. Moreover, we characterize the maximal elements of $\CP$.

If $\{Q_1,\dots,Q_n\}$ is the set formed by the maximal elements of $\mathcal C(\mathcal P)$, then a $\mathcal C(\mathcal P)$-set is a subset of $Q_i\setminus{\Delta(\mathcal C(\mathcal P))}$ for some $i\in\{1,\ldots,n\}$. In Section~\ref{sec:smallest} we prove that if $A$ is a $\mathcal C(\mathcal P)$-set, then there exists the smallest element of $\mathcal C(\mathcal P)$ that contains $A$. We denote this semigroup by $\mathcal C(\mathcal P)[A]$ and say that $\mathcal C(\mathcal P)[A]$ is the $\mathcal C(\mathcal P)$-semigroup generated by $A$. If $S=\mathcal C(\mathcal P)[A]$ and $S\ne\mathcal C(\mathcal P)[B]$ for all $B\subsetneq A$, then we say that $A$ is a $\mathcal C(\mathcal P)$-\textit{minimal system of generators} of $S$. In Section~\ref{sec:smallest} we show that every $\mathcal C(\mathcal P)$-semigroup admits a unique $\mathcal C(\mathcal P)$-minimal system of generators. If $S$ is a $\mathcal C(\mathcal P)$-semigroup, then we denote by $\mathcal C(\mathcal P)\mathrm{msg}(S)$ its $\mathcal C(\mathcal P)$-minimal system of generators. The cardinality of $\mathcal C(\mathcal P)\mathrm{msg}(S)$ is named the $\mathcal C(\mathcal P)$-rank of $S$. Finally, Section~\ref{sec:rankone} is devoted to the study of $\mathcal C(\mathcal P)$-semigroups of $\mathcal C(\mathcal P)$-rank 1.

We will see that the $\CP$-semigroups of $\CP$-rank one are exactly those $S$ for which $w-m(S)\in\langle\mathcal P\rangle$ for every element $w$ of $\Ap(S,m(S))$ with $0<w<F(S)$, and that all of them have maximal embedding dimension. As a consequence, a numerical semigroup $S$ with $m(S)<F(S)$ has maximal embedding dimension if and only if it has $\mathcal C(\PF(S))$-rank one.

\section{The Tree Associated with the Covariety \texorpdfstring{$\mathcal C(\mathcal P)$}{C(P)}}\label{sec:tree}

Our first objective is to prove that $\mathcal C(\mathcal P)$ is a covariety.

\begin{lemma}(\cite[Lemma 3.2]{covariedades})\label{lema primero}
    If $S$ and $T$ are numerical semigroups and $x\in S$, then:
    \begin{enumerate}
        \item $S\cap T$ is a numerical semigroup and $F(S\cap T)=\max\{F(S),F(T)\}$.
        \item $S\setminus\{x\}$ is a numerical semigroup if and only if $x\in\mathrm{msg}(S)$.
        \item $m(S)=\min(\mathrm{msg}(S))$.
    \end{enumerate}
\end{lemma}

\begin{lemma} \label{lema delta es el minimo}
    $\Delta=\{0,\max(\mathcal P)+1,\rightarrow\}$ is the minimum of $\mathcal C(\mathcal P)$.
\end{lemma}

\begin{proof}
    It is clear that $\Delta$ is a numerical semigroup, $F(\Delta)=\max(\mathcal P)$ and $\mathcal P\subseteq\mathrm{PF}(\Delta)$. Therefore, $\Delta\in\mathcal C(\mathcal P)$. If $S\in\mathcal C(\mathcal P)$, then $F(S)=\max(\mathcal P)$. Hence, $\Delta\subseteq S$.
\end{proof}

\begin{lemma}\label{lema interseccion pseudo-Frobenius}
    If $S$ and $T$ are numerical semigroups and $f\in\mathrm{PF}(S)\cap\mathrm{PF}(T)$, then $f\in\mathrm{PF}(S\cap T)$.
\end{lemma}

\begin{proof}
    Since $f\in\mathrm{PF}(S)$, then $f\notin S$ and therefore $f\notin S\cap T$. If $x\in (S\cap T)\setminus\{0\}$, then $f+x\in S$ and $f+x\in T$. Hence, $f+x\in S\cap T$. Consequently, $f\in\mathrm{PF}(S\cap T)$.
\end{proof}

\noindent The converse of Lemma \ref{lema interseccion pseudo-Frobenius} is false. Let $S=\langle 3,8,10\rangle=\{0,3,6,8,\rightarrow\}$ and $T=\langle 4,6,9,11\rangle=\{0,4,6,8,\rightarrow\}$. Then $\PF(S)=\{5,7\}$ and $\PF(T)=\{2,5,7\}$, while $S\cap T=\{0,6,8,\rightarrow\}$ and $\PF(S\cap T)=\{2,3,4,5,7\}$. As Lemma \ref{lema interseccion pseudo-Frobenius} guarantees, $5$ and $7$ are pseudo-Frobenius numbers of $S\cap T$. But so is $3$, which belongs to $S$ and is not a pseudo-Frobenius number of $T$, because $3+4=7\notin T$.

\begin{lemma}\label{lema pseudo-Frobenius sin multiplicidad}
    If $S$ is a numerical semigroup, $S\ne \N$ and $f\in\mathrm{PF}(S)$, then $f\in\mathrm{PF}(S\setminus\{m(S)\})$.
\end{lemma}

\begin{proof}
    By Lemma \ref{lema primero}, we know that $S\setminus\{m(S)\}$ is a numerical semigroup and, since $S\ne\N$, then $f\in\N\setminus\{0\}$. Since $f\notin S$, then $f\notin S\setminus\{m(S)\}$. If $s\in S\setminus\{0,m(S)\}$, then $f+s\in S$ and $f+s\neq m(S)$ because $s>m(S)$. Hence $f+s\in S\setminus\{m(S)\}$ and $f\in\mathrm{PF}(S\setminus\{m(S)\})$.
\end{proof}

We are now in a position to prove the result announced at the beginning of this section.

\begin{proposition}\label{proposicion C es una covariedad}
    $\mathcal C(\mathcal P)$ is a covariety and $\Delta(\mathcal C(\mathcal P))=\{0,\max(\mathcal P)+1,\rightarrow\}$.
\end{proposition}

\begin{proof}
    By Lemma \ref{lema delta es el minimo}, we know that $\{0,\max(\mathcal P)+1,\rightarrow\}$ is the minimum of $\mathcal C(\mathcal P)$. If $\{S,T\}\subseteq\mathcal C(\mathcal P)$, then $F(S)=F(T)=\max(\mathcal P)$ and $\mathcal P\subseteq\mathrm{PF}(S)\cap\mathrm{PF}(T)$. Applying Lemmas \ref{lema primero} and \ref{lema interseccion pseudo-Frobenius}, we obtain that $F(S\cap T)=\max(\mathcal P)$ and $\mathcal P\subseteq\mathrm{PF}(S\cap T)$. Therefore, $S\cap T\in\mathcal C(\mathcal P)$. If $S\in\mathcal C(\mathcal P)$ and $S\ne\Delta(\mathcal C(\mathcal P))$, then $m(S)<F(S)$ and so $F(S\setminus\{m(S)\})=F(S)=\max(\mathcal P)$. Moreover, by Lemma \ref{lema pseudo-Frobenius sin multiplicidad}, $\mathcal P\subseteq \mathrm{PF}(S\setminus\{m(S)\})$. Therefore, $S\setminus\{m(S)\}\in\mathcal C(\mathcal P)$.
\end{proof}

As a consequence of \cite[Proposition 2.1]{covariedades}, we obtain the following result.

\begin{corollary}
    $\mathcal C(\mathcal P)$ is a finite set.
\end{corollary}

\begin{remark}\label{nota frobenius}
If $F\in\N\setminus\{0\}$, then $\mathcal C(\{F\})$ is the set of all numerical semigroups with Frobenius number $F$, since $F(S)\in\PF(S)$ for every numerical semigroup $S\neq\N$. Hence $\mathcal C(\{F\})$ is the covariety studied in \cite{covariedades}, and the results of this paper extend those of \cite{covariedades}. Moreover, directly from the definition of $\CP$,
\[
\CP=\bigcap_{p\in\mathcal P}\mathcal C(\{p,\max(\mathcal P)\}),
\qquad\text{and}\qquad
\mathcal C(\mathcal Q)\subseteq\CP \text{ whenever } \mathcal P\subseteq\mathcal Q \text{ and } \max(\mathcal P)=\max(\mathcal Q).
\]
\end{remark}

A \textit{graph} $G$ is a pair $(\mathcal V,E)$ where $\mathcal{V}$ is a nonempty set and $E$ is a subset of $\{(u,v)\in\mathcal V\times\mathcal{V} \mid u\ne v\}$. The elements of $\mathcal V$ and $E$ are the \textit{vertices} and \textit{edges} of $G$, respectively. A \textit{path} (of length $n$) connecting the vertices $x$ and $y$ of $G$ is a sequence of distinct edges of $G$ of the form $(v_0,v_1),(v_1,v_2),\ldots,(v_{n-1},v_n)$ such that $v_0=x$ and $v_n=y$.

A graph $G$ is a \textit{tree} if there exists a vertex $r$ (known as the \textit{root} of $G$) such that for any other vertex $x$ of $G$ there exists a unique path connecting $x$ and $r$. If $(u,v)$ is an edge of a tree $G$, then we say that $u$ is a child of $v$. A \textit{leaf} is a vertex without children. If there is a path connecting $u$ and $v$, we say that $u$ is a \textit{descendant} of $v$. The \textit{level} (or depth) of a vertex $x$ is the length of the path connecting $x$ and $r$, and the \textit{height} of $G$ is the maximum of the levels of its vertices.

We define the graph $G(\mathcal C(\mathcal P))$ as follows: $\mathcal C(\mathcal P)$ is its set of vertices and $(S,T)\in\mathcal C(\mathcal P)\times\mathcal C(\mathcal P)$ is an edge if and only if $T=S\setminus\{m(S)\}$.

As a consequence of Proposition \ref{proposicion C es una covariedad} and \cite[Proposition 2.3]{covariedades}, we obtain the following result.

\begin{proposition}
    $G(\mathcal C(\mathcal P))$ is a tree with root $\Delta(\mathcal C(\mathcal P))$.
\end{proposition}

A tree can be constructed recursively starting from the root and connecting each vertex with its children by means of an edge. Therefore, it is of interest to characterize the children of an arbitrary vertex of the tree $G(\mathcal C(\mathcal P))$. To this end, we introduce some concepts and results.

Given a numerical semigroup $S$, we denote by $\mathrm{SG}(S)=\{x\in\mathrm{PF}(S) \mid 2x\in S\}$. The elements of $\mathrm{SG}(S)$ are the \textit{special gaps} of $S$.

\begin{lemma}(\cite[Proposition 4.33]{libro})
    Let $S$ be a numerical semigroup and $x\in\N\setminus S$. Then $x\in\mathrm{SG}(S)$ if and only if $S\cup\{x\}$ is a numerical semigroup.
\end{lemma}

As a consequence of Proposition \ref{proposicion C es una covariedad} and \cite[Proposition 2.4]{covariedades}, we obtain the following result.

\begin{proposition}\label{proposicion hijos arbol}
    If $S\in\mathcal C(\mathcal P)$, then the set formed by the children of $S$ in the tree $G(\mathcal C(\mathcal P))$ is $$\{S\cup\{x\} \mid x\in\mathrm{SG}(S), x<m(S), S\cup\{x\}\in\mathcal{ C(\mathcal P)}\}.$$
\end{proposition}

We now aim to determine the conditions that $x$ must satisfy in order to have $S\cup\{x\}\in\mathcal C(\mathcal P)$.

\begin{lemma}\label{lema union hueco especial}
    Let $S\in\mathcal C(\mathcal P)$ and $x\in\mathrm{SG}(S)$. Then $S\cup\{x\}\in\mathcal C(\mathcal P)$ if and only if $x\notin\mathcal P$ and $\mathcal P+\{x\}\subseteq S$.
\end{lemma}

\begin{proof}
    (Necessity) If $S\cup\{x\}\in\mathcal C(\mathcal P)$, then $\mathcal P\subseteq\mathrm{PF}(S\cup\{x\})$. Therefore, $x\notin\mathcal P$ and $\mathcal P+\{x\}\subseteq S\cup\{x\}$. Since $p+x\neq x$ for every $p\in\mathcal P$, we get $\mathcal P+\{x\}\subseteq S$.

\noindent (Sufficiency) Since $\max(\mathcal P)=F(S)\in\mathcal P$ and $x\notin\mathcal P$, we have $x<F(S)$ and so $F(S \cup \{x\}) = \max(\mathcal{P})$. Moreover, since $\mathcal{P} \subseteq \mathrm{PF}(S)$ and $x\notin\mathcal P$, then $\mathcal P\cap(S\cup\{x\})=\emptyset$. Let $p\in\mathcal P$ and $s\in (S\cup\{x\})\setminus\{0\}$. If $s\in S$, then $p+s\in S$ because $p\in\PF(S)$; if $s=x$, then $p+s\in S$ by hypothesis. Therefore, $\mathcal P\subseteq\mathrm{PF}(S\cup\{x\})$. Consequently, $S\cup \{x\}\in\mathcal C(\mathcal P)$.
\end{proof}

\begin{lemma}\label{lema suma de pseudo-Frobenius}
    Let $S$ be a numerical semigroup, $S\ne\N$ and $\{f,g\}\subseteq\mathrm{PF}(S)$. Then, $f+g\in S$ if and only if $f+g\notin\mathrm{PF}(S)$.
\end{lemma}

\begin{proof}
    If $s\in S\setminus\{0\}$, then $g+s\in S\setminus\{0\}$ and therefore $f+g+s\in S$. Consequently, $f+g\in S\cup\mathrm{PF}(S)$. Hence, since $S\cap\PF(S)=\emptyset$, $f+g\in S$ if and only if $f+g\notin\mathrm{PF}(S)$.
\end{proof}

As a consequence of Lemmas \ref{lema union hueco especial} and \ref{lema suma de pseudo-Frobenius}, we obtain the following result.

\begin{proposition}\label{proposicion interseccion P y pseudos}
    Let $S\in\mathcal C(\mathcal P)$ and $x\in\mathrm{SG}(S)$. Then $S\cup\{x\}\in\mathcal C(\mathcal P)$ if and only if $x\notin\mathcal P$ and $(\mathcal P+\{x\})\cap\mathrm{PF}(S)=\emptyset$.
\end{proposition}

As a consequence of Propositions \ref{proposicion hijos arbol} and \ref{proposicion interseccion P y pseudos}, we obtain the following result.

\begin{corollary}\label{corolario caracterizacion hijos}
    If $S\in\mathcal C(\mathcal P)$, then the set formed by the children of $S$ in the tree $G(\mathcal C(\mathcal P))$ is $$\{S\cup\{x\}\mid x\in(\mathrm{SG}(S)\setminus\mathcal P), x<m(S), (\mathcal P +\{x\})\cap \mathrm{PF}(S)=\emptyset\}.$$
\end{corollary}

\begin{example}\label{ejemplo arbol 57}
We compute $\mathcal C(\{5,7\})$ using Corollary \ref{corolario caracterizacion hijos}. A number appearing next to an edge $(S, T)$ means that  $T=S\setminus\{m(S)\}$ (equivalently, $S=T\cup\{m(S)\}$).
\[
\begin{tikzpicture}[
    >=Stealth,
    every node/.style={font=\small},
    node distance=1cm and 1cm
]
\node (top) {$\Delta=\{0,8,\rightarrow\}$};
\node (a) [below left=of top] {$\langle4,9,10,11\rangle$};
\node (b) [below right=of top] {$\langle6,8,9,10,11,13\rangle$};
\node (c) [below left=of b] {$\langle3,8,10\rangle$};
\node (d) [below right=of b] {$\langle4,6,9,11\rangle$};
\draw[->] (a) -- node[midway,left,fill=white,inner sep=1pt] {$4$} (top);
\draw[->] (b) -- node[midway,right,fill=white,inner sep=1pt] {$6$} (top);
\draw[->] (c) -- node[midway,left,fill=white,inner sep=1pt] {$3$} (b);
\draw[->] (d) -- node[midway,right,fill=white,inner sep=1pt] {$4$} (b);
\end{tikzpicture}
\]
\end{example}

The rest of this section collects two properties of $G(\CP)$ that we will use repeatedly. The first one says where a semigroup sits in the tree; the second one, that the set of pseudo-Frobenius numbers shrinks along every edge. For a numerical semigroup $S$ we write $n(S)=\mathrm{Card}(\{s\in S\mid s<F(S)\})$, so that $g(S)+n(S)=F(S)+1$.

\begin{proposition}\label{proposicion nivel}
Let $S\in\CP$. Then the level of $S$ in $G(\CP)$ is
\[
\max(\mathcal P)-g(S)=n(S)-1=\mathrm{Card}(\{s\in S\mid 0<s<\max(\mathcal P)\}).
\]
Moreover, if $s_1<s_2<\cdots<s_d$ are the elements of $S$ in the interval $(0,\max(\mathcal P))$, then the path connecting $S$ and $\Delta(\CP)$ is
\[
S\ \to\ S\setminus\{s_1\}\ \to\ S\setminus\{s_1,s_2\}\ \to\ \cdots\ \to\ S\setminus\{s_1,\ldots,s_d\}=\Delta(\CP).
\]
\end{proposition}

\begin{proof}
Every edge $(S,T)$ satisfies $T=S\setminus\{m(S)\}$, so $g(T)=g(S)+1$. Since $g(\Delta)=F(\Delta)=\max(\mathcal P)$, we deduce that the level of $S$ is $\max(\mathcal P)-g(S)$, and $\max(\mathcal P)-g(S)=n(S)-1$ because $F(S)=\max(\mathcal P)$. The description of the path follows from the fact that $s_1=m(S)$ and, inductively, $s_{i+1}=m(S\setminus\{s_1,\dots,s_i\})$.
\end{proof}

\begin{example}\label{ejemplo nivel}
The semigroup $S=\langle 5,7,8\rangle=\{0,5,7,8,10,12,\rightarrow\}$ has $\PF(S)=\{9,11\}$, so it belongs to $\mathcal C(\{9,11\})$. Its elements in $(0,11)$ are $5,7,8,10$. Therefore $S$ has level $4$ in $G(\mathcal C(\{9,11\}))$, its genus is $11-4=7$, and the path connecting $S$ and $\Delta=\{0,12,\rightarrow\}$ is
\[
\langle 5,7,8\rangle\to\langle 7,8,10,12,13\rangle\to\langle 8,10,12,13,14,15,17,19\rangle\to\langle 10,12,13,\ldots,19,21\rangle\to\Delta .
\]
\end{example}

\begin{lemma}\label{lema PF decrece}
Let $S\in\CP$ and let $T=S\cup\{x\}$ be a child of $S$ in $G(\CP)$. Then $\PF(T)\subseteq\PF(S)\setminus\{x\}$. In particular, $t(T)<t(S)$.
\end{lemma}

\begin{proof}
Let $f\in\PF(T)$. Then $f\notin T$, so $f\notin S$ and $f\ne x$. If $s\in S\setminus\{0\}$, then $f+s\in T$. By Corollary \ref{corolario caracterizacion hijos}, $x<m(S)\leq s<f+s$, so $f+s\neq x$ and therefore $f+s\in S$. Hence $f\in\PF(S)\setminus\{x\}$. Since $x\in\SG(S)\subseteq \PF(S)$, the inclusion $\PF(T)\subseteq\PF(S)\setminus\{x\}$ gives $t(T)\leq t(S)-1$.
\end{proof}

\begin{example}\label{ejemplo PF arista}
In $G(\mathcal C(\{9,11\}))$, the semigroup $\langle 8,10,12,13,14,15,17,19\rangle=\{0,8,10,12,\rightarrow\}$ has pseudo-Frobenius numbers $2,4,5,6,7,9,11$. Its child $\langle 7,8,10,12,13\rangle$, obtained by adding $7$, has pseudo-Frobenius numbers $5,6,9,11$. The integer $7$ is lost because it is now an element, and $2$ and $4$ are lost because $2+7=9$ and $4+7=11$ are not elements of the child.
\end{example}

\begin{remark}\label{nota hipotesis}
Lemma \ref{lema PF decrece} uses that the special gap added to $S$ is smaller than $m(S)$. If we add a special gap greater than $m(S)$, the inclusion may fail. For instance, $S=\{0,4,5,8,\rightarrow\}$ has $\PF(S)=\{6,7\}$, while $S\cup\{7\}=\langle 4,5,7\rangle$ has $\PF(S\cup\{7\})=\{3,6\}$, and $3\notin\PF(S)$.
\end{remark}

\begin{corollary}\label{corolario PF cadena}
If $S\in\CP$ and $T$ is a descendant of $S$ in $G(\CP)$, then $\mathcal P\subseteq\PF(T)\subsetneq\PF(S)$.
\end{corollary}

\section{The Algorithm to compute \texorpdfstring{$\mathcal C(\mathcal P)$}{C(P)}}\label{sec:algorithm}

In this section we explain how to obtain from $\Ap(S,n)$ the data that the construction of the tree requires, namely $m(S)$, $\PF(S)$, $\SG(S)$ and the Ap\'ery sets of the children of $S$. We start by recalling the precise description of $\Ap(S,n)$.

\begin{lemma}(\cite[Lemma 2.4]{libro})
    If $S$ is a numerical semigroup and $n\in S\setminus\{0\}$, then $\Ap(S,n)$ is a set of cardinality $n$. Moreover, $\mathrm{Ap}(S,n)=\{w(0),w(1),\ldots,w(n-1)\}$ where $w(i)$ is the smallest element of $S$ that is congruent to $i$ modulo $n$ for all $i\in\{0,\ldots,n-1\}$.
\end{lemma}

If $S$ is a numerical semigroup, then we define on $\mathbb Z$ the following order relation: $a\leq_S b$ if and only if $b-a\in S$.

\begin{lemma}(\cite[Lemma 10]{JPAA}, see also \cite[Proposition 2.20]{libro})\label{lema pseudo-Frobenius y maximales}
    If $S$ is a numerical semigroup and $n\in S\setminus\{0\}$, then $$\mathrm{PF}(S)=\{w-n \mid w\in\mathrm{Maximals}_{\leq_S}(\mathrm{Ap}(S,n))\}.$$
\end{lemma}

The following result has an immediate proof.

\begin{lemma}\label{lema apery con apery}
    Let $S$ be a numerical semigroup, $n\in S\setminus\{0\}$ and $w\in\mathrm{Ap}(S,n)$. Then
    \[
    w\in\mathrm{Maximals}_{\leq_S}(\mathrm{Ap}(S,n)) \quad\text{if and only if}\quad
    w+w'\notin \mathrm{Ap}(S,n)\quad\text{for all } w'\in\mathrm{Ap}(S,n)\setminus\{0\}.
    \]
\end{lemma}

The following result follows from Lemma \ref{lema suma de pseudo-Frobenius}.

\begin{lemma}\label{lema huecos especiales}
    If $S$ is a numerical semigroup and $S\ne\N$, then $\mathrm{SG}(S)=\{x\in\mathrm{PF}(S)\mid 2x\notin \mathrm{PF}(S)\}$.
\end{lemma}

\begin{remark}\label{nota primera}
    Observe that, as a consequence of Lemmas \ref{lema pseudo-Frobenius y maximales}, \ref{lema apery con apery} and \ref{lema huecos especiales}, if $S$ is a numerical semigroup and we know $\Ap(S,n)$ for some $n\in S\setminus\{0\}$, then we can easily compute $\mathrm{PF}(S)$ and $\mathrm{SG}(S)$. If moreover $n=F(S)+1$, then an integer $s$ with $0<s<n$ belongs to $S$ if and only if it belongs to $\Ap(S,n)$; in particular, $m(S)=\min\big((\Ap(S,n)\setminus\{0\})\cup\{n\}\big)$.
\end{remark}

We illustrate the content of Remark \ref{nota primera} with an example.

\begin{example}
    If $S=\{0,7,8,10,\rightarrow\}$, then $\mathrm{Ap}(S,7)=\{0,8,10,11,12,13,16\}$. Applying Lemma \ref{lema apery con apery}, we obtain that $\mathrm{Maximals}_{\leq_S}(\mathrm{Ap}(S,7))=\{10,11,12,13,16\}$. Then, Lemma \ref{lema pseudo-Frobenius y maximales} tells us that $\mathrm{PF}(S)=\{3,4,5,6,9\}$ and Lemma \ref{lema huecos especiales} that $\mathrm{SG}(S)=\{4,5,6,9\}$.
\end{example}

The following result has an immediate proof.

\begin{lemma}\label{lema apery y huecos especiales}
    If $S$ is a numerical semigroup, $n\in S\setminus\{0\}$ and $x\in\mathrm{SG}(S)$, then $x+n\in\mathrm{Ap}(S,n)$. Moreover, $\mathrm{Ap}(S\cup\{x\},n)=(\mathrm{Ap}(S,n)\setminus\{x+n\})\cup\{x\}$.
\end{lemma}

\begin{remark}\label{nota segunda}
    Observe that, as a consequence of Lemma \ref{lema apery y huecos especiales}, if we know $\mathrm{Ap}(S,n)$, then we can easily determine $\mathrm{Ap}(S\cup\{x\},n)$. In particular, applying Proposition \ref{proposicion hijos arbol}, we have that from $\mathrm{Ap}(S,n)$ we can easily compute $\mathrm{Ap}(T,n)$ for every child $T$ of $S$ in the tree $G(\mathcal C(\mathcal P))$.
\end{remark}

\begin{example}
   If $S=\{0,6,8,\rightarrow\}$, then, by applying Lemmas \ref{lema pseudo-Frobenius y maximales}, \ref{lema apery con apery} and \ref{lema huecos especiales}, we obtain $\mathrm{Ap}(S,6)=\{0,8,9,10,11,13\}$, $\mathrm{PF}(S)=\{2,3,4,5,7\}$ and $\mathrm{SG}(S)=\{3,4,5,7\}$. Now, using Lemma \ref{lema apery y huecos especiales}, we deduce that
$\mathrm{Ap}(S\cup\{3\},6)=\{0,3,8,10,11,13\}$, $\mathrm{Ap}(S\cup\{4\},6)=\{0,4,8,9,11,13\}$ and $\mathrm{Ap}(S\cup\{5\},6)=\{0,5,8,9,10,13\}$.
\end{example}

We are now in a position to present the algorithm that gives the title to this section. Since every $S\in\CP$ has $F(S)=\max(\mathcal P)$, the integer $n=\max(\mathcal P)+1$ belongs to every element of $\CP$, and we work with $\Ap(S,n)$ throughout. The algorithm computes the tree level by level: after $d$ iterations of the main loop, $B$ is the set of vertices of level $d$ of $G(\CP)$, which by Proposition \ref{proposicion nivel} are the elements of $\CP$ with genus $\max(\mathcal P)-d$.

\begin{algo}{Computing $\mathcal C(\mathcal P)$}{algoritmo primero}
\KwIn{A finite nonempty subset $\mathcal P$ of $\mathbb N\setminus\{0\}$.}
\KwOut{The set $\mathcal C(\mathcal P)$.}
\BlankLine
\Fn(\Comment*[f]{$W=\Ap(S,n)$ with $n=F(S)+1$; Remark \ref{nota primera}}){\FData{$W$}}{
  $m\gets\min\big((W\setminus\{0\})\cup\{n\}\big)$\;
  $\mathrm{Max}\gets\{w\in W\setminus\{0\}\mid w+w'\notin W \text{ for all } w'\in W\setminus\{0\}\}$\;
  $\mathrm{PF}\gets\{w-n\mid w\in \mathrm{Max}\}$\;
  $\mathrm{SG}\gets\{x\in\mathrm{PF}\mid 2x\notin \mathrm{PF}\}$\;
  \KwRet{$(m,\mathrm{PF},\mathrm{SG})$}\;
}
\BlankLine
$n\gets \max(\mathcal P)+1$, \quad $\Delta \gets \{0,n,\rightarrow\}$, \quad $\Ap(\Delta,n)\gets\{0,n+1,n+2,\ldots,2n-1\}$\;
$A\gets\{\Delta\}$, \quad $B\gets\{\Delta\}$\;
\Loop{
  \ForAll{$S\in B$}{
    $(m,\mathrm{PF},\mathrm{SG})\gets$ \FData{$\Ap(S,n)$}\;
    $\theta(S)\gets\{x\in\mathrm{SG}\setminus\mathcal P\mid x<m,\ (\mathcal P+\{x\})\cap\mathrm{PF}=\emptyset\}$\;
  }
  \If{$\bigcup_{S\in B}\theta(S)=\emptyset$}{
    \KwRet{$A$}\;
  }
  $C\gets\bigcup_{S\in B}\{S\cup\{x\}\mid x\in\theta(S)\}$\;
  \ForAll(\Comment*[f]{Remark \ref{nota segunda}}){$S\in B$ \KwSty{and} $x\in\theta(S)$}{
    $\Ap(S\cup\{x\},n)\gets(\Ap(S,n)\setminus\{x+n\})\cup\{x\}$\;
  }
  $A\gets A\cup C$, \quad $B\gets C$\;
}
\end{algo}

Since $G(\CP)$ is a tree, every element of $\CP$ is produced exactly once, so no comparisons between semigroups are needed. For $\mathcal P=\{F\}$, Algorithm \ref{algoritmo primero} computes all the numerical semigroups with Frobenius number $F$ (Remark \ref{nota frobenius}). We illustrate the functioning of Algorithm \ref{algoritmo primero} with an example.

\begin{example}\label{ejemplo segundo}
    We compute $\mathcal C(\{5,7\})$ using Algorithm \ref{algoritmo primero}. Here $n=8$ and $\Delta=\{0,8,\rightarrow\}$. The following table shows, for each semigroup produced, the data computed by the function \textup{\textsc{Data}} and the set $\theta(S)$. The column $d$ is the iteration in which $S$ belongs to $B$.
    \[
    \begin{array}{cllcccc}
    \toprule
    d & S & \Ap(S,8) & m(S) & \PF(S) & \SG(S) & \theta(S)\\
    \midrule
    0 & \Delta & \{0,9,10,11,12,13,14,15\} & 8 & \{1,\ldots,7\} & \{4,5,6,7\} & \{4,6\}\\
    1 & \Delta\cup\{4\} & \{0,4,9,10,11,13,14,15\} & 4 & \{5,6,7\} & \{5,6,7\} & \emptyset\\
    1 & \Delta\cup\{6\} & \{0,6,9,10,11,12,13,15\} & 6 & \{2,3,4,5,7\} & \{3,4,5,7\} & \{3,4\}\\
    2 & \Delta\cup\{3,6\} & \{0,3,6,9,10,12,13,15\} & 3 & \{5,7\} & \{5,7\} & \emptyset\\
    2 & \Delta\cup\{4,6\} & \{0,4,6,9,10,11,13,15\} & 4 & \{2,5,7\} & \{2,5,7\} & \emptyset\\
    \bottomrule
    \end{array}
    \]
    For $\Delta\cup\{4,6\}$ the special gap $2$ is smaller than the multiplicity, but $2+5=7\in\PF(\Delta\cup\{4,6\})$, so it is discarded. In the first iteration $B=\{\Delta\}$ and we obtain $C=\{\Delta\cup\{4\},\Delta\cup\{6\}\}$; in the second one $B=C$ and we obtain $C=\{\Delta\cup\{3,6\},\Delta\cup\{4,6\}\}$; in the third one $\theta(\Delta\cup\{3,6\})=\theta(\Delta\cup\{4,6\})=\emptyset$, so the algorithm returns
    \[
    \mathcal C(\{5,7\})=\{\Delta,\Delta\cup\{4\},\Delta\cup\{6\},\Delta\cup\{3,6\},\Delta\cup\{4,6\}\}.
    \]
\end{example}

\section{\texorpdfstring{Some applications of $G(\mathcal{C}(\mathcal{P}))$}{Some applications of G(C(P))}}\label{sec:applications}

In this section we address some of the questions arising both from the construction of $G(\mathcal{C}(\mathcal{P}))$ and from Algorithm \ref{algoritmo primero}. Throughout, we write $F=\max(\mathcal P)$.

The first observation is as follows:

\begin{remark}\label{nota PF igual P}
    Observe that Algorithm \ref{algoritmo primero} allows us to compute the set
\[
\{S \mid S \text{ is a numerical semigroup and } \mathrm{PF}(S)=\mathcal{P}\}:
\]
it consists of the elements of $\CP$ of type $\mathrm{Card}(\mathcal P)$. In particular, the algorithm decides whether this set is empty. This problem was solved by Delgado, Garc\'{\i}a-S\'anchez and Robles-P\'erez \cite{delgado}, whose algorithm is implemented in the function \NSPF{} of the package \texttt{numericalsgps} \cite{numericalsgps-pkg,numericalsgps}. In Subsection~\ref{subsec:PF} we will see that, in order to compute this set with the tree, it is not necessary to build the whole of $G(\CP)$: the results on the type obtained in Subsection~\ref{subsec:type} allow us to discard whole subtrees.
\end{remark}

\begin{example}
    In Example \ref{ejemplo segundo} we obtained that $\mathrm{PF}(\Delta)=\{1,2,3,4,5,6,7\}$, $\mathrm{PF}(\Delta\cup\{4\})=\{5,6,7\}$, $\mathrm{PF}(\Delta\cup\{6\})=\{2,3,4,5,7\}$, $\mathrm{PF}(\Delta\cup\{3,6\})=\{5,7\}$ and $\PF(\Delta\cup\{4,6\})=\{2,5,7\}$. Consequently, we can conclude that $\{S\mid S\text{ is a numerical semigroup and } \PF(S)=\{5,7\}\}=\{\Delta\cup\{3,6\}\}=\{\langle 3,8,10\rangle\}$.
\end{example}

\subsection{\texorpdfstring{Semigroups with a prescribed type}{Semigroups with a prescribed type}}\label{subsec:type}

A natural question, which is difficult to answer \emph{a priori}, is the following: does there exist a numerical semigroup having a prescribed set $\mathcal{P}$ of pseudo-Frobenius numbers and a fixed type? The tree associated with the covariety $\mathcal{C}(\mathcal{P})$ provides an answer to this question. For this purpose, let us consider a remark, together with a variation of Algorithm \ref{algoritmo primero} and some examples:

\begin{remark}\label{nota tipo decrese}
    In $G(\mathcal{C}(\mathcal{P}))$, the type strictly decreases along each branch, by Lemma \ref{lema PF decrece} (see also Remark \ref{nota hipotesis}).
\end{remark}

Before giving the algorithm, we locate the semigroups of a given type inside the tree. We need two classical inequalities. The first one is due to Fr\"oberg, Gottlieb and H\"aggkvist \cite{froberg}: for every numerical semigroup $S\neq\N$,
\begin{equation}\label{eq FGH}
F(S)+1\leq (t(S)+1)\,n(S).
\end{equation}
The second one is due to Nari \cite[Proposition 2.2]{nari}: for every numerical semigroup $S$,
\begin{equation}\label{eq Nari}
2g(S)\geq F(S)+t(S),
\end{equation}
and, by \cite[Proposition-Definition 2.3]{nari}, equality holds if and only if $S$ is \textit{almost symmetric} in the sense of Barucci and Fr\"oberg \cite{barucci-froberg}.

\begin{proposition}\label{proposicion banda}
If $S\in\CP$ has type $k$ and level $d$ in $G(\CP)$, then
\[
\left\lceil\frac{F+1}{k+1}\right\rceil-1\ \leq\ d\ \leq\ \left\lfloor\frac{F-k}{2}\right\rfloor .
\]
\end{proposition}

\begin{proof}
By Proposition \ref{proposicion nivel}, $d=F-g(S)=n(S)-1$. The lower bound follows from \eqref{eq FGH} and the upper bound from \eqref{eq Nari}.
\end{proof}

\begin{corollary}\label{corolario altura}
The height of $G(\CP)$ is at most $\left\lfloor\frac{F-\mathrm{Card}(\mathcal P)}{2}\right\rfloor$, and every $S\in\CP$ satisfies $g(S)\geq \left\lceil\frac{F+\mathrm{Card}(\mathcal P)}{2}\right\rceil$. Moreover, there exists $S\in\CP$ with $2g(S)=F+\mathrm{Card}(\mathcal P)$ if and only if there exists an almost symmetric numerical semigroup $S$ with $\PF(S)=\mathcal P$.
\end{corollary}

\begin{proof}
The first two statements follow from Proposition \ref{proposicion banda}, since $t(S)\geq \mathrm{Card}(\mathcal P)$ for all $S\in\CP$. If $2g(S)=F+\mathrm{Card}(\mathcal P)$, then by \eqref{eq Nari} we have $F+t(S)\leq F+\mathrm{Card}(\mathcal P)$, so $t(S)=\mathrm{Card}(\mathcal P)$, hence $\PF(S)=\mathcal P$, and equality holds in \eqref{eq Nari}, so $S$ is almost symmetric. The converse is clear.
\end{proof}

The bound of Corollary \ref{corolario altura} is attained in Examples \ref{ejemplo arbol 57} and \ref{ejemplo 9 11}, but not in general. For instance, $\mathcal C(\{4,6\})=\{\{0,7,\rightarrow\},\langle 5,7,8,9,11\rangle\}$ has height $1<\left\lfloor\frac{6-2}{2}\right\rfloor$.

\begin{example}\label{ejemplo 9 11}
Let $\mathcal P=\{9,11\}$, so $F=11$, $n=12$ and $\Delta=\{0,12,\rightarrow\}$. Table \ref{tabla 9 11} shows the whole tree $G(\mathcal C(\{9,11\}))$, computed with Algorithm \ref{algoritmo primero}. For each vertex $S$ it gives its level $d$, its parent $S\setminus\{m(S)\}$ and its pseudo-Frobenius numbers; by Proposition \ref{proposicion nivel}, the genus of $S$ is $11-d$. It has $23$ vertices and height $4=\lfloor (11-2)/2\rfloor$. In agreement with Proposition \ref{proposicion banda}, the semigroups of type $3$ appear in levels $2$, $3$ and $4$, and those of type $2$ only in level $4$.
\end{example}

{\small
\setlength{\tabcolsep}{5pt}
\renewcommand{\arraystretch}{0.9}
\setlength{\LTpre}{\smallskipamount}
\begin{longtable}{@{}cclccl@{}}
\caption{The tree $G(\mathcal C(\{9,11\}))$.}\label{tabla 9 11}\\
\toprule
 & $d$ & $\msg(S)$ & parent & $t(S)$ & $\PF(S)$\\
\midrule
\endfirsthead
\multicolumn{6}{@{}l}{\emph{Table \ref{tabla 9 11} (continued)}}\\
\toprule
 & $d$ & $\msg(S)$ & parent & $t(S)$ & $\PF(S)$\\
\midrule
\endhead
\bottomrule
\endlastfoot
$S_{0}$ & 0 & $\langle 12, 13, \ldots, 23\rangle$ & -- & 11 & $\{1,2,\ldots,11\}$\\
$S_{1}$ & 1 & $\langle 6, 13, 14, 15, 16, 17\rangle$ & $S_{0}$ & 5 & $\{7, 8, 9, 10, 11\}$\\
$S_{2}$ & 1 & $\langle 7, 12, 13, 15, 16, 17, 18\rangle$ & $S_{0}$ & 6 & $\{5, 6, 8, 9, 10, 11\}$\\
$S_{3}$ & 1 & $\langle 8, 12, 13, 14, 15, 17, 18, 19\rangle$ & $S_{0}$ & 7 & $\{4, 5, 6, 7, 9, 10, 11\}$\\
$S_{4}$ & 1 & $\langle 10, 12, 13, \ldots, 19, 21\rangle$ & $S_{0}$ & 9 & $\{2, 3, \ldots, 9, 11\}$\\
$S_{5}$ & 2 & $\langle 4, 13, 14, 15\rangle$ & $S_{3}$ & 3 & $\{9, 10, 11\}$\\
$S_{6}$ & 2 & $\langle 5, 12, 13, 14, 16\rangle$ & $S_{4}$ & 4 & $\{7, 8, 9, 11\}$\\
$S_{7}$ & 2 & $\langle 6, 7, 15, 16, 17\rangle$ & $S_{2}$ & 4 & $\{8, 9, 10, 11\}$\\
$S_{8}$ & 2 & $\langle 6, 8, 13, 15, 17\rangle$ & $S_{3}$ & 4 & $\{7, 9, 10, 11\}$\\
$S_{9}$ & 2 & $\langle 6, 10, 13, 14, 15, 17\rangle$ & $S_{4}$ & 5 & $\{4, 7, 8, 9, 11\}$\\
$S_{10}$ & 2 & $\langle 7, 8, 12, 13, 17, 18\rangle$ & $S_{3}$ & 5 & $\{5, 6, 9, 10, 11\}$\\
$S_{11}$ & 2 & $\langle 7, 10, 12, 13, 15, 16, 18\rangle$ & $S_{4}$ & 6 & $\{3, 5, 6, 8, 9, 11\}$\\
$S_{12}$ & 2 & $\langle 8, 10, 12, 13, 14, 15, 17, 19\rangle$ & $S_{4}$ & 7 & $\{2, 4, 5, 6, 7, 9, 11\}$\\
$S_{13}$ & 3 & $\langle 4, 10, 13, 15\rangle$ & $S_{12}$ & 3 & $\{6, 9, 11\}$\\
$S_{14}$ & 3 & $\langle 5, 7, 13, 16\rangle$ & $S_{11}$ & 3 & $\{8, 9, 11\}$\\
$S_{15}$ & 3 & $\langle 5, 8, 12, 14\rangle$ & $S_{12}$ & 3 & $\{7, 9, 11\}$\\
$S_{16}$ & 3 & $\langle 6, 7, 8, 17\rangle$ & $S_{10}$ & 3 & $\{9, 10, 11\}$\\
$S_{17}$ & 3 & $\langle 6, 7, 10, 15\rangle$ & $S_{11}$ & 3 & $\{8, 9, 11\}$\\
$S_{18}$ & 3 & $\langle 6, 8, 10, 13, 15, 17\rangle$ & $S_{12}$ & 5 & $\{2, 4, 7, 9, 11\}$\\
$S_{19}$ & 3 & $\langle 7, 8, 10, 12, 13\rangle$ & $S_{12}$ & 4 & $\{5, 6, 9, 11\}$\\
$S_{20}$ & 4 & $\langle 4, 6, 13, 15\rangle$ & $S_{18}$ & 3 & $\{2, 9, 11\}$\\
$S_{21}$ & 4 & $\langle 5, 7, 8\rangle$ & $S_{19}$ & 2 & $\{9, 11\}$\\
$S_{22}$ & 4 & $\langle 6, 7, 8, 10\rangle$ & $S_{19}$ & 2 & $\{9, 11\}$\\
\end{longtable}
}

The band of Proposition \ref{proposicion banda} restricts the levels that must be explored. The next result restricts the vertices: for each $S\in\CP$ it gives pseudo-Frobenius numbers of $S$ that remain pseudo-Frobenius numbers in all the descendants of $S$. Recall from Corollary \ref{corolario caracterizacion hijos} that the children of $S$ are obtained by adding to $S$ an integer smaller than $m(S)$.

\begin{definition}\label{definicion admisible}
Let $S\in\CP$. We say that an integer $x$ is \textit{admissible} for $S$ if $2\leq x<m(S)$, $x\notin\mathcal P$ and, for every $p\in\mathcal P$, we have $p-x\notin S$, and $p+x\in S$ whenever $p+x>m(S)$. We denote by $\PFs(S)$ the set of pseudo-Frobenius numbers $f$ of $S$ that are not admissible for $S$ and satisfy $f+x\in S$ for every integer $x$ admissible for $S$.
\end{definition}

Both notions depend on $\mathcal P$, which is fixed. By Remark \ref{nota primera}, the integers admissible for $S$ and the set $\PFs(S)$ are easily computed from $\Ap(S,\max(\mathcal P)+1)$.

\begin{lemma}\label{lema congelados}
Let $S\in\CP$ and let $T$ be a descendant of $S$ in $G(\CP)$, or $T=S$. Then:
\begin{enumerate}
\item every element of $T\setminus S$ is admissible for $S$. In particular, $S$ and $T$ have the same elements greater than $m(S)$.
\item $\mathcal P\cup\PFs(S)\subseteq\PF(T)$.
\end{enumerate}
\end{lemma}

\begin{proof}
(1) We may assume $T\neq S$. Let $S=T_0, T_1,\ldots,T_r=T$ be the vertices of the path connecting $T$ and $S$, so that $T_i=T_{i-1}\cup\{x_i\}$ for some $x_i<m(T_{i-1})$ (Corollary \ref{corolario caracterizacion hijos}). Since $m(T_i)=x_i$, we get $x_r<\cdots<x_1<m(S)$, and so $T\setminus S=\{x_1,\ldots,x_r\}\subseteq(0,m(S))$. Moreover, $1\notin T$ because $T\neq\N$. Let $x\in T\setminus S$. Then $x\notin\mathcal P$, because $\mathcal P\cap T=\emptyset$. Let $p\in\mathcal P$. If $p-x\in S$, then $p-x\neq0$ and $p=(p-x)+x\in T$, which is impossible. Since $p\in\PF(T)$, we have $p+x\in T$; if $p+x>m(S)$, then $p+x\in S$, because $T\setminus S\subseteq(0,m(S))$. Hence $x$ is admissible for $S$.

(2) We know that $\mathcal P\subseteq\PF(T)$, because $T\in\CP$. Let $f\in\PFs(S)$. Since $f\notin S$ and $f$ is not admissible for $S$, (1) gives $f\notin T$. If $t\in T\setminus\{0\}$, then either $t\in S$, and $f+t\in S$ because $f\in\PF(S)$; or $t\in T\setminus S$, so that $t$ is admissible for $S$ by (1), and $f+t\in S$ by the definition of $\PFs(S)$. In both cases $f+t\in T$. Hence $f\in\PF(T)$.
\end{proof}

\begin{corollary}\label{corolario poda}
Let $S\in\CP$ and $k\in\N$. If $\mathrm{Card}(\mathcal P\cup\PFs(S))>k$, then neither $S$ nor any of its descendants in $G(\CP)$ has type $k$. In particular, if $\PFs(S)\not\subseteq\mathcal P$, then no descendant $T$ of $S$ satisfies $\PF(T)=\mathcal P$.
\end{corollary}

\begin{example}\label{ejemplo poda}
Let $\mathcal P=\{9,11\}$ and $S=S_3=\{0,8,12,\rightarrow\}$ (Table \ref{tabla 9 11}). Then $m(S)=8$ and $\PF(S)=\{4,5,6,7,9,10,11\}$. The integer $2$ is not admissible for $S$, because $9+2=11\notin S$, and neither is $3$, because $11-3=8\in S$; on the other hand, $4$, $5$, $6$ and $7$ are admissible for $S$. Hence $\PFs(S)=\{9,10,11\}$. By Lemma \ref{lema congelados}, $10$ is a pseudo-Frobenius number of every descendant of $S_3$, which agrees with Table \ref{tabla 9 11}: the descendants of $S_3$ are $S_5$, $S_8$, $S_{10}$ and $S_{16}$. Therefore, if we look for the semigroups $T$ with $\PF(T)=\{9,11\}$, the subtree rooted at $S_3$ can be discarded.
\end{example}

We now rewrite the algorithm that computes the semigroups of $\CP$ with a fixed type. By Remark \ref{nota tipo decrese}, it is not necessary to explore the descendants of a semigroup of type smaller than or equal to $k$; by Proposition \ref{proposicion banda}, it is enough to explore the levels $0,1,\ldots,\lfloor (F-k)/2\rfloor$; and by Corollary \ref{corolario poda}, the descendants of $S$ can be discarded when $\mathrm{Card}(\mathcal P\cup\PFs(S))>k$.

\begin{algo}{Computing the semigroups in $G(\mathcal C(\mathcal P))$ with fixed type}{algoritmo tipo}
\KwIn{A finite nonempty subset $\mathcal P$ of $\mathbb N\setminus\{0\}$ and a natural number $k\ne0$.}
\KwOut{The set $\{S\in\mathcal C(\mathcal P)\mid t(S)=k\}$.}
$F\gets\max(\mathcal P)$, \quad $n\gets F+1$\;
\If{$k<\mathrm{Card}(\mathcal P)$ \KwSty{or} $k>F$}{
  \KwRet{$\emptyset$}\;
}
$D\gets\lfloor (F-k)/2\rfloor$\Comment*[r]{last level that can contain type $k$}
$B\gets\{\Delta\}$ with $\Delta=\{0,n,\rightarrow\}$, \quad $R\gets\emptyset$, \quad $d\gets 0$\;
\While{$B\neq\emptyset$ \KwSty{and} $d\leq D$}{
  $C\gets\emptyset$\;
  \ForAll{$S\in B$}{
    $(m,\mathrm{PF},\mathrm{SG})\gets$ \FData{$\Ap(S,n)$}\Comment*[r]{as in Algorithm \ref{algoritmo primero}}
    \uIf{$\mathrm{Card}(\mathrm{PF})=k$}{
      $R\gets R\cup\{S\}$\Comment*[r]{its descendants have type $<k$}
    }
    \ElseIf{$\mathrm{Card}(\mathrm{PF})>k$ \KwSty{and} $d<D$ \KwSty{and} $\mathrm{Card}(\mathcal P\cup\PFs(S))\leq k$}{
      $\theta(S)\gets\{x\in\mathrm{SG}\setminus\mathcal P\mid x<m,\ (\mathcal P+\{x\})\cap\mathrm{PF}=\emptyset\}$\;
      $C\gets C\cup\{S\cup\{x\}\mid x\in\theta(S)\}$\Comment*[r]{update Ap as in Algorithm \ref{algoritmo primero}}
    }
  }
  $B\gets C$, \quad $d\gets d+1$\;
}
\KwRet{$R$}\;
\end{algo}

\begin{example}\label{ejemplo tipo 4}
  We continue with Example \ref{ejemplo segundo}. Suppose that we want to determine whether there exists a numerical semigroup $S$ such that $\{5,7\}\subseteq\mathrm{PF}(S)$, $F(S)=7$ and $t(S)=4$. Algorithm \ref{algoritmo tipo} only needs to explore the levels $0$ and $1$, because $\lfloor (7-4)/2\rfloor=1$. It expands the root, of type $7$, and computes its two children, of types $3$ and $5$. Neither child is expanded: the first one has type smaller than $4$, and the second one lies in level $1=D$. So no such semigroup exists. The figure shows the whole tree $G(\mathcal C(\{5,7\}))$ with the type of each vertex. The algorithm builds the three vertices drawn in black; the two vertices in gray are never computed.
\[
\begin{tikzpicture}[
    >=Stealth,
    every node/.style={font=\small, text=black},
    node distance=1cm and 1cm,
    nb/.style={text=black!40}
]
\node (top) {$\Delta=\{0,8,\rightarrow\}$};
\node[above=1mm of top] {$t=7$};
\node (a) [below left=of top] {$\langle4,9,10,11\rangle$};
\node[above=1mm of a] {$t=3$};
\node (b) [below right=of top] {$\langle6,8,9,10,11,13\rangle$};
\node[above=1mm of b] {$t=5$};
\node[nb] (c) [below left=of b] {$\langle3,8,10\rangle$};
\node[nb,above=1mm of c] {$t=2$};
\node[nb] (d) [below right=of b] {$\langle4,6,9,11\rangle$};
\node[nb,above=1mm of d] {$t=3$};
\draw[->,black] (a) -- (top);
\draw[->,black] (b) -- (top);
\draw[->,black!30] (c) -- (b);
\draw[->,black!30] (d) -- (b);
\end{tikzpicture}
\]
\end{example}

\begin{example}\label{ejemplo tipo 7}
Let $\mathcal P=\{9,11\}$ (Example \ref{ejemplo 9 11}). Suppose that we want to know whether there exists a numerical semigroup $S$ with $\{9,11\}\subseteq\PF(S)$, $F(S)=11$ and $t(S)=7$. By Proposition \ref{proposicion banda}, such a semigroup lies in level $1$ or $2$, so Algorithm \ref{algoritmo tipo} stops after level $2$. In level $1$ it finds $S_3$ (type $7$), discards $S_1$ and $S_2$ (types $5$ and $6$) and expands only $S_4$ (type $9$). Among the children of $S_4$ it finds $S_{12}$. Hence
\[
\{S\in\mathcal C(\{9,11\})\mid t(S)=7\}=\{\langle 8,12,13,14,15,17,18,19\rangle,\ \langle 8,10,12,13,14,15,17,19\rangle\}.
\]
The algorithm only builds $9$ of the $23$ vertices of the tree: $S_0,\ldots,S_4$ and the four children of $S_4$. The figure shows the whole tree $G(\mathcal C(\{9,11\}))$ with the type of each vertex (see Table \ref{tabla 9 11}). The vertices built by the algorithm are drawn in black, and the two semigroups that it returns are framed. The dashed line separates the levels $d\leq D=2$, the only ones that the algorithm explores.
\[
\begin{tikzpicture}[
    >=Stealth, x=1.15cm, y=1.75cm,
    every node/.style={font=\small, text=black},
    v/.style={align=center, inner sep=1.6pt},
    bt/.style={v},
    rt/.style={v, draw, rounded corners=2pt, thick},
    nb/.style={v, text=black!35},
    eb/.style={->, black},
    en/.style={->, black!22}
]
\node[bt] (s0) at (3.81,0) {$S_{0}$\\[-2pt]{\scriptsize $t=11$}};
\node[bt] (s1) at (0,-1) {$S_{1}$\\[-2pt]{\scriptsize $t=5$}};
\node[bt] (s2) at (1,-1) {$S_{2}$\\[-2pt]{\scriptsize $t=6$}};
\node[rt] (s3) at (3,-1) {$S_{3}$\\[-2pt]{\scriptsize $t=7$}};
\node[bt] (s4) at (7.625,-1) {$S_{4}$\\[-2pt]{\scriptsize $t=9$}};
\node[nb] (s7) at (1,-2) {$S_{7}$\\[-2pt]{\scriptsize $t=4$}};
\node[nb] (s5) at (2,-2) {$S_{5}$\\[-2pt]{\scriptsize $t=3$}};
\node[nb] (s8) at (3,-2) {$S_{8}$\\[-2pt]{\scriptsize $t=4$}};
\node[nb] (s10) at (4,-2) {$S_{10}$\\[-2pt]{\scriptsize $t=5$}};
\node[bt] (s6) at (5,-2) {$S_{6}$\\[-2pt]{\scriptsize $t=4$}};
\node[bt] (s9) at (6,-2) {$S_{9}$\\[-2pt]{\scriptsize $t=5$}};
\node[bt] (s11) at (7,-2) {$S_{11}$\\[-2pt]{\scriptsize $t=6$}};
\node[rt] (s12) at (10.25,-2) {$S_{12}$\\[-2pt]{\scriptsize $t=7$}};
\node[nb] (s16) at (4,-3) {$S_{16}$\\[-2pt]{\scriptsize $t=3$}};
\node[nb] (s14) at (6.5,-3) {$S_{14}$\\[-2pt]{\scriptsize $t=3$}};
\node[nb] (s17) at (7.5,-3) {$S_{17}$\\[-2pt]{\scriptsize $t=3$}};
\node[nb] (s13) at (8.5,-3) {$S_{13}$\\[-2pt]{\scriptsize $t=3$}};
\node[nb] (s15) at (9.5,-3) {$S_{15}$\\[-2pt]{\scriptsize $t=3$}};
\node[nb] (s18) at (10.5,-3) {$S_{18}$\\[-2pt]{\scriptsize $t=5$}};
\node[nb] (s19) at (12,-3) {$S_{19}$\\[-2pt]{\scriptsize $t=4$}};
\node[nb] (s20) at (10.5,-4) {$S_{20}$\\[-2pt]{\scriptsize $t=3$}};
\node[nb] (s21) at (11.5,-4) {$S_{21}$\\[-2pt]{\scriptsize $t=2$}};
\node[nb] (s22) at (12.5,-4) {$S_{22}$\\[-2pt]{\scriptsize $t=2$}};
\foreach \c/\p in {s1/s0,s2/s0,s3/s0,s4/s0,s6/s4,s9/s4,s11/s4,s12/s4}
  \draw[eb] (\c) -- (\p);
\foreach \c/\p in {s5/s3,s7/s2,s8/s3,s10/s3,s13/s12,s14/s11,s15/s12,s16/s10,s17/s11,s18/s12,s19/s12,s20/s18,s21/s19,s22/s19}
  \draw[en] (\c) -- (\p);
\draw[densely dashed, black!50] (-0.5,-2.5) -- (13,-2.5);
\node[font=\scriptsize, text=black!60, anchor=south west] at (-0.5,-2.5) {$D=2$};
\end{tikzpicture}
\]
\end{example}

An interesting question arising from Remark \ref{nota tipo decrese} is the determination of the minimal type among the semigroups in the tree $G(\mathcal C(\mathcal P))$. The following result provides a first step in this direction.

\begin{corollary}\label{corolario tipo minimo}
    We have
$$
\min\{t(S)\mid S\in \mathcal C(\mathcal P)\}
=
\min\{t(S)\mid S \text{ is a leaf in } G(\mathcal C(\mathcal P))\}.
$$
\end{corollary}

\begin{proof}
Every $S\in\CP$ has a leaf among its descendants (or is a leaf), and by Lemma \ref{lema PF decrece} the type of that leaf is at most $t(S)$.
\end{proof}

When there exists a numerical semigroup $S$ in the tree $G(\mathcal{C}(\mathcal{P}))$ such that $\mathrm{PF}(S)=\mathcal{P}$, it is clear that
\[
\mathrm{Minimals}\{\mathrm{PF}(S)\mid S\in\mathcal{C}(\mathcal{P})\}
=
\min\{\mathrm{PF}(S)\mid S\in\mathcal{C}(\mathcal{P})\}
=
\mathcal{P}.
\]
However, when no semigroup satisfies this condition, we obtain the following:

\begin{corollary}
    We have
$$\mathrm{Minimals}\{\mathrm{PF}(S)\mid S\in\mathcal C(\mathcal P)\}=\mathrm{Minimals}\{\mathrm{PF}(S) \mid S \text{ is a leaf in } G(\mathcal C(\mathcal P))\}.$$
\end{corollary}

\begin{proof}
It follows from Corollary \ref{corolario PF cadena}, with the same argument as in the proof of Corollary \ref{corolario tipo minimo}.
\end{proof}

\subsection{Semigroups with a prescribed set of pseudo-Frobenius numbers}\label{subsec:PF}

We now come back to the problem mentioned in Remark \ref{nota PF igual P}: to compute the numerical semigroups $S$ with $\PF(S)=\mathcal P$, which is what the function \NSPF{} does. The results of Subsection~\ref{subsec:type} solve it on a pruned version of $G(\CP)$.

\begin{proposition}\label{proposicion PF igual P}
Let $S$ be a numerical semigroup. The following conditions are equivalent:
\begin{enumerate}
\item $\PF(S)=\mathcal P$.
\item $S\in\CP$ and $t(S)=\mathrm{Card}(\mathcal P)$.
\end{enumerate}
In this case $S$ is a leaf of $G(\CP)$, and its level $d$ satisfies
\[
\left\lceil\frac{F+1}{\mathrm{Card}(\mathcal P)+1}\right\rceil-1\ \leq\ d\ \leq\ \left\lfloor\frac{F-\mathrm{Card}(\mathcal P)}{2}\right\rfloor .
\]
\end{proposition}

\begin{proof}
The equivalence is clear, because $\mathcal P\subseteq \PF(S)$ for every $S\in\CP$. If $\PF(S)=\mathcal P$, then $\SG(S)\setminus\mathcal P=\emptyset$, so $S$ has no children by Corollary \ref{corolario caracterizacion hijos}. The bounds are Proposition \ref{proposicion banda} with $k=\mathrm{Card}(\mathcal P)$.
\end{proof}

Therefore, Algorithm \ref{algoritmo tipo} with input $\mathcal P$ and $k=\mathrm{Card}(\mathcal P)$ returns
\[
\{S\mid S \text{ is a numerical semigroup and } \PF(S)=\mathcal P\}.
\]
This is the same set that \NSPF{} returns, but Algorithm \ref{algoritmo tipo} obtains it without building the whole tree. If the output is empty, then no numerical semigroup has $\mathcal P$ as its set of pseudo-Frobenius numbers. Note that for $k=\mathrm{Card}(\mathcal P)$ the bound $D=\lfloor (F-k)/2\rfloor$ coincides with the bound on the height of $G(\CP)$ given in Corollary \ref{corolario altura}, so it does not discard any level; in this case the reduction of the search comes from Corollary \ref{corolario poda}, which discards the subtree rooted at $S$ as soon as $\PFs(S)\not\subseteq\mathcal P$.

\begin{example}\label{ejemplo 4 9}
Is there a numerical semigroup $S$ with $\PF(S)=\{4,9\}$? Here $F=9$, $\mathrm{Card}(\mathcal P)=2$, and Proposition \ref{proposicion PF igual P} says that such an $S$ would be in level $3$ of $G(\mathcal C(\{4,9\}))$, since $\lceil 10/3\rceil-1=3=\lfloor 7/2\rfloor$. The whole tree is the following (the type is written above each semigroup; the two vertices in gray are not built by Algorithm \ref{algoritmo tipo}, as explained below):
\[
\begin{tikzpicture}[
    >=Stealth,
    every node/.style={font=\small, text=black},
    nb/.style={text=black!40}
]
\node[align=center] (r) at (0,0) {$t=9$\\ $\langle 10,11,\ldots,19\rangle$};
\node[align=center] (a) at (-5,-1.9) {$t=5$\\ $\langle 6,10,11,13,14,15\rangle$};
\node[align=center] (b) at (0,-1.9) {$t=6$\\ $\langle 7,10,11,12,13,15,16\rangle$};
\node[align=center] (c) at (5,-1.9) {$t=7$\\ $\langle 8,10,11,12,13,14,15,17\rangle$};
\node[align=center,nb] (d) at (0,-3.8) {$t=4$\\ $\langle 6,7,10,11,15\rangle$};
\node[align=center] (e) at (3.5,-3.8) {$t=5$\\ $\langle 6,8,10,11,13,15\rangle$};
\node[align=center] (f) at (7.5,-3.8) {$t=5$\\ $\langle 7,8,10,11,12,13\rangle$};
\node[align=center,nb] (g) at (7.5,-5.7) {$t=3$\\ $\langle 6,7,8,10,11\rangle$};
\draw[->,black] (a)--(r); \draw[->,black] (b)--(r); \draw[->,black] (c)--(r);
\draw[->,black!30] (d)--(b); \draw[->,black] (e)--(c); \draw[->,black] (f)--(c);
\draw[->,black!30] (g)--(f);
\end{tikzpicture}
\]
The only semigroup in level $3$ is $\langle 6,7,8,10,11\rangle$, whose set of pseudo-Frobenius numbers is $\{4,5,9\}$. Hence there is no numerical semigroup with $\PF(S)=\{4,9\}$, and the minimum type in $\mathcal C(\{4,9\})$ is $3$. Corollary \ref{corolario poda} shortens the search. For $S=\langle 7,10,11,12,13,15,16\rangle=\{0,7,10,\rightarrow\}$ we have $\PF(S)=\{3,4,5,6,8,9\}$, the integers admissible for $S$ are $3$ and $6$, and $\PFs(S)=\{4,8,9\}$; since $8\notin\{4,9\}$, the subtree rooted at $S$ is discarded, and the same happens with $\langle 7,8,10,11,12,13\rangle$, whose only child is $\langle 6,7,8,10,11\rangle$. Thus Algorithm \ref{algoritmo tipo} visits $6$ of the $8$ vertices of the tree (those drawn in black) and returns $\emptyset$.

For $\mathcal P=\{9,11\}$ the band of Proposition \ref{proposicion PF igual P} is $\{3,4\}$, and Table \ref{tabla 9 11} shows that
\[
\{S\mid \PF(S)=\{9,11\}\}=\{\langle 5,7,8\rangle,\ \langle 6,7,8,10\rangle\}.
\]
With the pruning of Corollary \ref{corolario poda} (see Example \ref{ejemplo poda}), only $17$ of the $23$ vertices of $G(\mathcal C(\{9,11\}))$ are visited by Algorithm \ref{algoritmo tipo}.
\end{example}

\subsection{Maximal elements of \texorpdfstring{$\CP$}{C(P)}}\label{subsec:maximal}

The maximal elements of $\CP$ appear in the definition of $\CP$-set, so it is useful to identify them. Every maximal element is a leaf of $G(\CP)$, but the converse is false.

\begin{proposition}\label{proposicion maximales}
Let $S\in\CP$. Then $S$ is a maximal element of $\CP$ (with respect to inclusion) if and only if there is no $x\in\SG(S)\setminus\mathcal P$ such that $(\mathcal P+\{x\})\cap\PF(S)=\emptyset$. In particular, if $\PF(S)=\mathcal P$, then $S$ is a maximal element of $\CP$.
\end{proposition}

\begin{proof}
If such an $x$ exists, then $S\subsetneq S\cup\{x\}\in\CP$ by Proposition \ref{proposicion interseccion P y pseudos}. Conversely, assume that $S\subsetneq T$ for some $T\in\CP$, and let $x=\max(T\setminus S)$. If $s\in S\setminus\{0\}$, then $x+s\in T$ and $x+s>x$, so $x+s\in S$; hence $x\in\PF(S)$. Analogously $2x\in T$ and $2x>x$, so $2x\in S$ and $x\in\SG(S)$. Since $\mathcal P\cap T=\emptyset$, we have $x\notin\mathcal P$; and $\mathcal P+\{x\}\subseteq T$ with $p+x>x$ for every $p\in\mathcal P$, so $\mathcal P+\{x\}\subseteq S$. As $S\cap\PF(S)=\emptyset$, we conclude that $(\mathcal P+\{x\})\cap\PF(S)=\emptyset$.
\end{proof}

\begin{example}
In $\mathcal C(\{5,7\})$ the maximal elements are $\langle 3,8,10\rangle$ and $\langle 4,6,9,11\rangle$. The leaf $\langle 4,9,10,11\rangle$ is not maximal, since $6\in\SG(\langle 4,9,10,11\rangle)$ and $\langle 4,9,10,11\rangle\cup\{6\}=\langle 4,6,9,11\rangle$. Note also that the maximal element $\langle 4,6,9,11\rangle$ has $\PF=\{2,5,7\}\neq\{5,7\}$. In $\mathcal C(\{9,11\})$ there are $14$ leaves but only $3$ maximal elements, $S_{20}$, $S_{21}$ and $S_{22}$ (Table \ref{tabla 9 11}). Finally, $\mathcal C(\{4,9\})$ has a maximum, $\langle 6,7,8,10,11\rangle$, so every subset of $\{6,7,8\}$ is a $\mathcal C(\{4,9\})$-set.
\end{example}

\subsection{\texorpdfstring{Semigroups with a prescribed genus}{Semigroups with a prescribed genus}}\label{subsec:genus}

A further question of interest is the determination, within the tree $G(\mathcal{C}(\mathcal{P}))$, of all numerical semigroups with a prescribed genus. From the construction of $G(\mathcal{C}(\mathcal{P}))$, it follows that the genus decreases by one unit along each edge, from the root to the leaves. By Proposition \ref{proposicion nivel}, the semigroups of genus $k$ are exactly those of level $F-k$, and by Corollary \ref{corolario altura} there are none if $k<\lceil (F+\mathrm{Card}(\mathcal P))/2\rceil$. This observation leads to the following procedure, which provides yet another modification of Algorithm \ref{algoritmo primero}.

\begin{algo}{Computing the semigroups with fixed genus in $G(\mathcal C(\mathcal P))$}{algoritmo genero}
\KwIn{A finite nonempty subset $\mathcal P$ of $\mathbb N\setminus\{0\}$ and a nonzero natural number $k$.}
\KwOut{$\{S\in\mathcal C(\mathcal P)\mid g(S)=k\}$.}
$F\gets\max(\mathcal P)$\;
\If{$k>F$ \KwSty{or} $k<\lceil (F+\mathrm{Card}(\mathcal P))/2\rceil$}{
  \KwRet{$\emptyset$}\;
}
Run Algorithm \ref{algoritmo primero} during $F-k$ iterations of its main loop\;
\KwRet{the set $B$ obtained after these iterations, or $\emptyset$ if Algorithm \ref{algoritmo primero} stops earlier}\;
\end{algo}

\begin{example}
For $\mathcal P=\{9,11\}$ and $k=8$, Algorithm \ref{algoritmo genero} returns level $3$ of Table \ref{tabla 9 11}, that is, the seven semigroups $S_{13},\ldots,S_{19}$:
\begin{gather*}
\langle 4,10,13,15\rangle,\ \langle 5,7,13,16\rangle,\ \langle 5,8,12,14\rangle,\ \langle 6,7,8,17\rangle,\\ \langle 6,7,10,15\rangle,\ \langle 6,8,10,13,15,17\rangle,\ \langle 7,8,10,12,13\rangle.
\end{gather*}
For $k=6$ the algorithm returns $\emptyset$ without computing anything, since $6<\lceil 13/2\rceil=7$.
\end{example}

\begin{remark}
The minimum genus in $\CP$ is $F-h$, where $h$ is the height of $G(\CP)$, and the semigroups attaining it are the elements of the set $B$ at the moment in which Algorithm \ref{algoritmo primero} returns $A$. In $\mathcal C(\{9,11\})$ the minimum genus is $7$, attained by $S_{20}$, $S_{21}$ and $S_{22}$; here the bound of Corollary \ref{corolario altura} is attained, but there is no almost symmetric semigroup with $\PF(S)=\{9,11\}$ because $11+2$ is odd.
\end{remark}

\section{The smallest \texorpdfstring{$\mathcal C(\mathcal P)$}{C(P)}-semigroup containing a \texorpdfstring{$\mathcal C(\mathcal P)$}{C(P)}-set}\label{sec:smallest}

Since $\CP$ is finite, every $\CP$-semigroup is contained in a maximal element of $\CP$ with respect to inclusion; these maximal elements are characterized in Proposition \ref{proposicion maximales}. As in the Introduction, let $Q_1,\ldots,Q_n$ be the maximal elements of $\CP$. A \textit{$\CP$-set} is a subset of $Q_i\setminus\Delta(\CP)$ for some $i\in\{1,\ldots,n\}$. In other words, a set $A$ of integers is a $\CP$-set if and only if it is contained in some $\CP$-semigroup and all its elements are positive integers smaller than $\max(\mathcal P)$. For instance, the $\mathcal C(\{9,11\})$-sets are the subsets of $S_{20}\setminus\Delta=\{4,6,8,10\}$, of $S_{21}\setminus\Delta=\{5,7,8,10\}$ and of $S_{22}\setminus\Delta=\{6,7,8,10\}$ (see Table \ref{tabla 9 11}).

If $A$ is a $\mathcal C(\mathcal P)$-set, then we denote by $\mathcal C(\mathcal P)[A]$ the intersection of all $\mathcal{ C(\mathcal P)}$-semigroups containing $A$. Since $\mathcal C(\mathcal P)$ is a finite set, it is closed under intersection, and $A$ is contained in at least one maximal element of $\CP$, the intersection of the $\mathcal C(\mathcal P)$-semigroups containing $A$ is a $\mathcal C(\mathcal P)$-semigroup. Therefore, we can state the following result.

\begin{proposition}
    If $A$ is a $\mathcal C(\mathcal P)$-set, then $\mathcal C(\mathcal P)[A]$ is the smallest $\mathcal C(\mathcal P)$-semigroup containing $A$.
\end{proposition}

The following result follows from \cite[Proposition 4.2]{covariedades}.

\begin{proposition}\label{proposicion covariedad conjunto}
    If $S\in\mathcal C(\mathcal P)$, then $A=\{x\in\mathrm{msg}(S)\mid x\notin\Delta(\mathcal C(\mathcal P))\}$ is a $\mathcal C(\mathcal P)$-set and $\mathcal C(\mathcal P)[A]=S$.
\end{proposition}

If $A$ is a $\mathcal C(\mathcal P)$-set and $S=\mathcal C(\mathcal P)[A]$, then we say that $A$ is a $\mathcal C(\mathcal P)$-system of generators of $S$. If, in addition, $S\ne\mathcal C(\mathcal P)[B]$ for all $B\subsetneq A$, then we say that $A$ is a $\mathcal C(\mathcal P)$-minimal system of generators of $S$.

\begin{example}\label{ejemplo generadores 57}
    In Example \ref{ejemplo segundo} we saw that $\mathcal C(\{5,7\})=\{\Delta=\langle 8,9,10,11,12,13,14,15\rangle, \Delta\cup\{4\}=\langle 4,9,10,11\rangle, \Delta\cup\{6\}=\langle6,8,9,10,11,13\rangle,\Delta\cup\{3,6\}=\langle 3,8,10 \rangle,\Delta\cup\{4,6\}=\langle 4,6,9,11 \rangle\}$. Using Proposition \ref{proposicion covariedad conjunto}, we have that $\Delta=\mathcal C(\{5,7\})[\emptyset]$, $\Delta\cup\{4\}$$=\mathcal C(\{5,7\})[\{4\}]$, $\Delta\cup\{6\}=\mathcal C(\{5,7\})[\{6\}]$, $\Delta\cup\{3,6\}=\mathcal C(\{5,7\})[\{3\}]$ and $\Delta\cup\{4,6\}=\mathcal C(\{5,7\})[\{4,6\}]$.
\end{example}

In Example \ref{ejemplo generadores 57}, the set $\{x\in\mathrm{msg}(S)\mid x\notin\Delta\}$ is a $\mathcal C(\{5,7\})$-minimal system of generators of $S$ for every $S\in\mathcal C(\{5,7\})$. This is not true in general, as Example \ref{ejemplo 3 5 17} below shows.

The following result has an immediate proof.

\begin{lemma}\label{lema covariedad de conjuntos inclusion}
    If $\mathcal X$ and $\mathcal Y$ are $\mathcal C(\mathcal P)$-sets and $\mathcal X\subseteq\mathcal Y$, then $\mathcal C(\mathcal P)[\mathcal X]\subseteq\mathcal C(\mathcal P)[\mathcal Y]$.
\end{lemma}

\begin{example}\label{ejemplo 3 5 17}
    Let $S=\langle 10,13,15,16,18,19,21,22,24,27 \rangle\in\mathcal C(\{3,5,17\})$. By Proposition \ref{proposicion covariedad conjunto}, we know that $S=\mathcal C(\{3,5,17\})[\{10,13,15,16\}]$. However, $\{10,13,15,16\}$ is not a $\mathcal C(\{3,5,17\})$-minimal system of generators, since $S=\mathcal C(\{3,5,17\})[\{10\}]$. Indeed, since $\{10\}\subseteq \{10,13,15,16\}$, then by Lemma \ref{lema covariedad de conjuntos inclusion}, we have that $\mathcal C(\{3,5,17\})[\{10\}]\subseteq S$. Now, if $T$ is a $\mathcal C(\{3,5,17\})$-semigroup and $10\in T$, then $10+3=13$, $10+5=15$, $13+3=16$ and every integer greater than $17$ belongs to $T$. Therefore, $S\subseteq\mathcal C(\{3,5,17\})[\{10\}]$.
\end{example}

Our next objective in this section is to prove that every $\mathcal C(\mathcal P)$-semigroup admits a unique $\mathcal C(\mathcal P)$-minimal system of generators. To this end, we need to introduce some results.

The following result follows from \cite[Lemma 5.1]{covariedades}.

\begin{lemma}
    If $S\in\mathcal C(\mathcal P)$, $S\ne\Delta(\mathcal C(\mathcal P))$ and $A$ is a $\mathcal C(\mathcal P)$-set such that $\mathcal C(\mathcal P)[A]=S$, then $m(S)\in A$.
\end{lemma}

\begin{lemma}\label{lema covariedad x+s}
    If $S$ is a $\mathcal C(\mathcal P)$-semigroup, $s\in S\setminus\{0\}$ and $x\in\langle\mathcal P\rangle$, then $s+x\in S$.
\end{lemma}

\begin{proof}
    If $x=0$, the result is trivial. If $x\ne0$, then there exist $n\in\N\setminus\{0\}$ and $x_1,\ldots,x_n\in\mathcal P$ such that $x=x_1+\cdots+x_n$. We prove by induction on $n$ that $s+x\in S$. If $n=1$, the result holds since $\mathcal P\subseteq \mathrm{PF}(S)$. By the induction hypothesis, $s+x_1+\cdots+x_{n-1}\in S$. Since $x_n\in\mathrm{PF}(S)$, we have that $x_n+(s+x_1+\cdots+x_{n-1})\in S$ and, therefore, $s+x\in S$.
\end{proof}

We define on $\mathbb Z$ the following order relation: $a\leq_\mathcal P b$ if and only if $b-a\in\langle\mathcal P\rangle$.

\begin{lemma}\label{lema minimales y B}
    If $S\in\mathcal C(\mathcal P)$ and $B=\mathrm{Minimals}_{\leq_{\mathcal P}}\{x\in\mathrm{msg}(S) \mid x\notin\Delta(\mathcal C(\mathcal P))\}$, then $\mathcal C(\mathcal P)[B]=S$.
\end{lemma}

\begin{proof}
    It is clear that $B$ is a $\mathcal C(\mathcal P)$-set. Since $B\subseteq\{x\in\mathrm{msg}(S)\mid x\notin\Delta(\mathcal C(\mathcal P)) \}$, applying Lemma \ref{lema covariedad de conjuntos inclusion} and Proposition \ref{proposicion covariedad conjunto}, we obtain that $\mathcal C(\mathcal P)[B]\subseteq \mathcal C(\mathcal P)[\{x\in\mathrm{msg}(S)\mid x\notin\Delta(\mathcal C(\mathcal P))\}]=S$. Let us now prove the reverse inclusion. If $x\in\mathrm{msg}(S)$ is such that $x\notin\Delta(\mathcal C(\mathcal P))$, then there exist $b\in B$ and $c\in\langle\mathcal P\rangle$ such that $x=b+c$. Applying Lemma \ref{lema covariedad x+s}, we deduce that if $T\in\mathcal C(\mathcal P)$ and $B\subseteq T$, then $\{x\in\mathrm{msg}(S)\mid x\notin\Delta(\mathcal C(\mathcal P))\}\subseteq T$. Therefore, $S\subseteq\mathcal C(\mathcal P)[B]$.
\end{proof}

The following result is straightforward to prove.

\begin{lemma}
    Let $A$ be a $\mathcal C(\mathcal P)$-set and $S=\mathcal C(\mathcal P)[A]$. Then $A$ is a $\mathcal C(\mathcal P)$-minimal system of generators of $S$ if and only if $a\notin\mathcal C(\mathcal P)[A\setminus\{a\}]$ for all $a\in A$.
\end{lemma}

\begin{lemma}\label{lema minimales y B sin b}
    Let $S\in\mathcal C(\mathcal P)$ and $B=\mathrm{Minimals}_{\leq_\mathcal P}\{x\in\mathrm{msg}(S)\mid x\notin\Delta(\mathcal C(\mathcal P))\}$. If $b\in B$, then $S\setminus\{b\}\in\mathcal C(\mathcal P)$.
\end{lemma}

\begin{proof}
    Since $b\in\mathrm{msg}(S)$ and $b<F(S)=\max(\mathcal P)$, applying Lemma \ref{lema primero}, we deduce that $S\setminus\{b\}$ is a numerical semigroup and $F(S\setminus\{b\})=\max(\mathcal P)$. To conclude the proof, let us show that $\mathcal P+(S\setminus\{0,b\})\subseteq S\setminus\{b\}$. It suffices to prove that if $s\in S\setminus\{0\}$ and $p\in\mathcal P$, then $p+s\ne b$. We distinguish two cases:
    \begin{enumerate}
        \item If $s\notin\mathrm{msg}(S)$, then $s=s_1+s_2$ for some $s_1,s_2\in S\setminus\{0\}$. If $b=p+s$, then $b=(p+s_1)+s_2$ with $p+s_1\in S\setminus\{0\}$ and $s_2\in S\setminus\{0\}$, which contradicts the fact that $b\in\mathrm{msg}(S)$.
        \item If $s\in\mathrm{msg}(S)$ and $b=p+s$, then we conclude that $b\notin\mathrm{Minimals}_{\leq_\mathcal P}\{x\in\mathrm{msg}(S)\mid x\notin \Delta(\mathcal C(\mathcal P))\}$ since $s\in\mathrm{msg}(S)$, $s\notin\Delta(\mathcal C(\mathcal P))$ (because $s<b<\max(\mathcal P)$) and $s\leq_{\mathcal P}b$.
    \end{enumerate}
\end{proof}

We are now in a position to prove the result announced above.

\begin{theorem}\label{teorema}
    If $S\in \mathcal C(\mathcal P)$, then $B=\mathrm{Minimals}_{\leq_\mathcal P}\{x\in\mathrm{msg}(S)\mid x\notin\Delta(\mathcal C(\mathcal P))\}$ is the unique $\mathcal C(\mathcal P)$-minimal system of generators of $S$.
\end{theorem}

\begin{proof}
    By Lemma \ref{lema minimales y B}, we know that $\mathcal C(\mathcal P)[B]=S$. To conclude the proof, let us show that if $C$ is a $\mathcal C(\mathcal P)$-set and $\mathcal C(\mathcal P)[C]=S$, then $B\subseteq C$. If $B\not\subseteq C$, then there exists $b\in B\setminus C$. By Lemma \ref{lema minimales y B sin b}, we know that $S\setminus\{b\}\in\mathcal C(\mathcal P)$. Since $C\subseteq S\setminus\{b\}$, it follows that $\mathcal C(\mathcal P)[C]\ne S$, which is a contradiction.
\end{proof}

If $S\in\mathcal C(\mathcal P)$, then by Theorem \ref{teorema} we know that $S$ admits a unique $\mathcal C(\mathcal P)$-minimal system of generators, which we denote by $\mathcal C(\mathcal P)\mathrm{msg}(S)$. The cardinality of $\CP\msg(S)$ is called the $\mathcal C(\mathcal P)$-rank of $S$ and is denoted by $\CP\rk(S)$.

\begin{example}\label{ejemplo msg}
In Example \ref{ejemplo 3 5 17}, Theorem \ref{teorema} gives directly
\[
\mathcal C(\{3,5,17\})\msg(S)=\mathrm{Minimals}_{\leq_{\mathcal P}}\{10,13,15,16\}=\{10\},
\]
because $13=10+3$, $15=10+5$ and $16=10+3+3$. In $\mathcal C(\{9,11\})$ the situation is the opposite: two integers $x>y$ with $2\leq y<x<11$ satisfy $x-y\leq 8$, so $x-y\notin\langle 9,11\rangle$, and therefore $\mathcal C(\{9,11\})\msg(S)=\{x\in\msg(S)\mid x<11\}$ for every $S\in\mathcal C(\{9,11\})$ (see Table \ref{tabla 9 11}).
\end{example}

\begin{corollary}\label{corolario quitar}
    Let $S\in\mathcal C(\mathcal P)$ and $x\in S$. Then $S\setminus\{x\}\in\mathcal C(\mathcal P)$ if and only if $x\in \mathcal C(\mathcal P)\mathrm{msg}(S)$.
\end{corollary}

\begin{proof}
    (Necessity) Assume that $S\setminus\{x\}\in\CP$ and $x\notin\CP\msg(S)$. Then $\CP\msg(S)\subseteq S\setminus\{x\}$, and therefore $S=\CP[\CP\msg(S)]\subseteq S\setminus\{x\}$, which is a contradiction.

    (Sufficiency) It is a consequence of Theorem \ref{teorema} and Lemma \ref{lema minimales y B sin b}.
\end{proof}

Corollary \ref{corolario quitar} has a nice interpretation in terms of the partial order $(\CP,\subseteq)$. Recall that $S$ \textit{covers} $T$ if $T\subsetneq S$ and there is no $U\in\CP$ with $T\subsetneq U\subsetneq S$.

\begin{proposition}\label{proposicion Hasse}
Let $S\in\CP$. The elements of $\CP$ covered by $S$ are exactly the semigroups $S\setminus\{x\}$ with $x\in\CP\msg(S)$. In particular, $\CP\rk(S)$ is the number of elements of $\CP$ covered by $S$.
\end{proposition}

\begin{proof}
By Corollary \ref{corolario quitar}, each $S\setminus\{x\}$ with $x\in\CP\msg(S)$ is in $\CP$, and it is obviously covered by $S$. Conversely, let $T\in\CP$ with $T\subsetneq S$, and let $x=\min(S\setminus T)$. If $x=a+b$ with $a,b\in S\setminus\{0\}$, then $a,b<x$, so $a,b\in T$ and $x\in T$, a contradiction; hence $x\in\msg(S)$, and $T\subseteq S\setminus\{x\}$. Moreover, $x<F$ because $\Delta(\CP)\subseteq T$. If $p+s=x$ with $p\in\mathcal P$ and $s\in S\setminus\{0\}$, then $s<x$, so $s\in T$ and $x=p+s\in T$, a contradiction. Therefore $\mathcal P+(S\setminus\{0,x\})\subseteq S\setminus\{x\}$, which shows that $S\setminus\{x\}\in\CP$, and by Corollary \ref{corolario quitar}, $x\in\CP\msg(S)$. If $T$ is covered by $S$, then $T=S\setminus\{x\}$.
\end{proof}

So the tree $G(\CP)$ is a spanning tree of the Hasse diagram of $(\CP,\subseteq)$: each $S\neq\Delta(\CP)$ covers $\CP\rk(S)$ elements, and the tree keeps only one of these edges, the one corresponding to $m(S)$. For example, $S_{19}=\langle 7,8,10,12,13\rangle$ in Table \ref{tabla 9 11} has $\mathcal C(\{9,11\})$-rank $3$, and it covers $S_{19}\setminus\{7\}=S_{12}$ (its parent), $S_{19}\setminus\{8\}=S_{11}$ and $S_{19}\setminus\{10\}=S_{10}$.

\section{The \texorpdfstring{$\mathcal C(\mathcal P)$}{C(P)}-\texorpdfstring{semigroups}{semigroups} of \texorpdfstring{$\mathcal C(\mathcal P)$}{C(P)}-rank one}\label{sec:rankone}

Our objective in this section is to characterize the $\mathcal C(\mathcal P)$-semigroups of $\mathcal C(\mathcal P)$-rank one.

The following result is a consequence of \cite[Propositions 6.1 and 6.2]{covariedades}.

\begin{lemma}
    If $S\in \mathcal C(\mathcal P)$, then the following conditions hold:
    \begin{enumerate}
        \item $\CP\rk(S)\leq e(S)$.
        \item $\CP\rk(S)=0$ if and only if $S=\Delta(\mathcal C(\mathcal P))$.
    \end{enumerate}
\end{lemma}

The following result is a consequence of \cite[Proposition 6.3]{covariedades}.

\begin{lemma}\label{lema caracterizancion rango 1}
    If $S\in\mathcal C(\mathcal P)$ and $S\ne\Delta(\mathcal C(\mathcal P))$, then the following conditions are equivalent:
    \begin{enumerate}
        \item $\CP\rk(S)=1$.
        \item There exists $a\in S$ such that $\mathcal C(\mathcal P)[\{a\}]=S$.
        \item $S=\mathcal C(\mathcal P)[\{m(S)\}]$.
    \end{enumerate}
\end{lemma}

\begin{lemma}\label{lema M es submonoide}
    If $m\in \N\setminus\{0\}$, then $M=(\{m\}+\langle\mathcal P\cup\{m\}\rangle)\cup\{0\}$ is a submonoid of $(\N,+)$.
\end{lemma}

\begin{proof}
    It is clear that $M=\{km+x\mid k\in\N\setminus\{0\}, x\in\langle\mathcal P\rangle\}\cup\{0\}$. If $\{k_1,k_2\}\subseteq\N\setminus\{0\}$ and $\{x_1,x_2\}\subseteq\langle\mathcal P\rangle$, then
    \[
    k_1m+x_1+k_2m+x_2=(k_1+k_2)m+(x_1+x_2),
    \]
    where $k_1+k_2\in\mathbb N\setminus\{0\}$ and $x_1+x_2\in\langle\mathcal P\rangle$. Therefore, $M$ is closed under addition and contains zero. Hence, $M$ is a submonoid of $(\N,+)$.
\end{proof}

\begin{lemma}\label{lema M contenido en S}
    If $S$ is a $\mathcal C(\mathcal P)$-semigroup, $m\in S\setminus\{0\}$ and $M=(\{m\}+\langle\mathcal P\cup\{m\}\rangle)\cup\{0\}$, then $M\subseteq S$.
\end{lemma}

\begin{proof}
    If $a\in M\setminus\{0\}$, then $a=km+x$ with $k\in\N\setminus\{0\}$ and $x\in\langle\mathcal P\rangle$. Since $km\in S\setminus\{0\}$, by Lemma \ref{lema covariedad x+s}, we obtain that $a\in S$.
\end{proof}

\begin{lemma}\label{lema caracterizacion T}
    Let $m\in\N\setminus\{0,1\}$. Then
    \[
    T=\{0\}\cup(\{m\}+\langle\mathcal P\cup\{m\}\rangle)\cup\{\max(\mathcal P)+1,\rightarrow\}
    \]
    is a $\mathcal{C}(\mathcal P)$-semigroup if and only if $\mathcal P\cap T=\emptyset$.
\end{lemma}

\begin{proof}
    By Lemma \ref{lema M es submonoide}, we deduce that $T$ is a numerical semigroup.

    (Necessity) If $T\in\mathcal C(\mathcal P)$, then $\mathcal P\subseteq\mathrm{PF}(T)$ and therefore $\mathcal P\cap T=\emptyset$.

    (Sufficiency) If $\mathcal P\cap T=\emptyset$, then $F(T)=\max(\mathcal P)$. To conclude the proof, it suffices to show that $\mathcal P + (T\setminus\{0\})\subseteq T$. Let $p\in\mathcal P$ and $t\in T\setminus\{0\}$. We distinguish two cases:
    \begin{enumerate}
        \item If $t\geq \max(\mathcal P)+1$, then clearly $p+t\in T$.
        \item If $t<\max(\mathcal P)+1$, then $t\in\{m\}+\langle\mathcal P\cup\{m\}\rangle$ and clearly $p+t\in\{m\}+\langle \mathcal P\cup\{m\} \rangle\subseteq T$.
    \end{enumerate}
\end{proof}

\begin{proposition}\label{proposcion T rango uno}
    If $m\in\N\setminus\{0,1\}$, $m<\max(\mathcal P)$ and $(\{m\}+\langle\mathcal P\cup\{m\}\rangle)\cap\mathcal P=\emptyset$, then
    \[
    T=\{0\}\cup(\{m\}+\langle\mathcal P\cup\{m\}\rangle)\cup\{\max(\mathcal P)+1,\rightarrow\}
    \]
    is a $\mathcal C(\mathcal P)$-semigroup of $\mathcal C(\mathcal P)$-rank one. Moreover, every $\mathcal C(\mathcal P)$-semigroup of $\mathcal C(\mathcal P)$-rank one is of this form.
\end{proposition}

\begin{proof}
    Since $\mathcal P\cap T=\emptyset$, Lemma \ref{lema caracterizacion T} tells us that $T$ is a $\CP$-semigroup. As $m<\max(\mathcal P)$, we have $m\in T\setminus\Delta(\CP)$, so $\{m\}$ is a $\CP$-set. From Lemma \ref{lema M contenido en S}, we deduce that $T$ is the smallest $\mathcal C(\mathcal P)$-semigroup containing $\{m\}$. Therefore, $T=\mathcal C(\mathcal P)[\{m\}]$ and consequently $\CP\rk(T)=1$.

    If $S$ is a $\mathcal C(\mathcal P)$-semigroup and $\CP\rk(S)=1$, then by Lemma \ref{lema caracterizancion rango 1} we know that $S=\mathcal C(\mathcal P)[\{m(S)\}]$. Applying Lemmas \ref{lema M contenido en S} and \ref{lema caracterizacion T}, we deduce that
    \[
    S=\{0\}\cup(\{m(S)\}+\langle\mathcal P\cup\{m(S)\}\rangle)\cup\{\max(\mathcal P)+1,\rightarrow\},
    \]
    with $1<m(S)<\max(\mathcal P)$ and $(\{m(S)\}+\langle\mathcal P\cup\{m(S)\}\rangle)\cap \mathcal P=\emptyset$.
\end{proof}

\subsection{Rank one and maximal embedding dimension}\label{subsec:MED}

In this subsection we show that the semigroups of $\CP$-rank one have maximal embedding dimension, and that every MED-semigroup is of rank one in the covariety associated with its own set of pseudo-Frobenius numbers. We first recall a well-known characterization of MED-semigroups (see \cite[Chapter 2]{libro} and \cite{med}).

\begin{lemma}\label{lema MED tipo}
Let $S\neq\N$ be a numerical semigroup with multiplicity $m$. The following conditions are equivalent:
\begin{enumerate}
\item $S$ has maximal embedding dimension.
\item $t(S)=m-1$.
\item $\PF(S)=\{w-m\mid w\in\Ap(S,m)\setminus\{0\}\}$.
\end{enumerate}
\end{lemma}

\begin{lemma}\label{lema sumas de PF}
Let $S\neq\N$ be a numerical semigroup. If $x\in\langle\PF(S)\rangle$, $x\neq0$ and $x\notin S$, then $x\in\PF(S)$.
\end{lemma}

\begin{proof}
Write $x=f_1+\cdots+f_k$ with $f_i\in\PF(S)$; we use induction on $k$. If $k=1$ there is nothing to prove. If $k\geq 2$, let $y=f_1+\cdots+f_{k-1}$. If $y\in S$, then $y\neq0$ and $x=y+f_k\in S$, a contradiction. Hence $y\notin S$ and, by the induction hypothesis, $y\in\PF(S)$. By the proof of Lemma \ref{lema suma de pseudo-Frobenius}, $x=y+f_k\in S\cup\PF(S)$, and so $x\in\PF(S)$.
\end{proof}

\begin{theorem}\label{teorema MED}
Let $S\in\CP$ with $S\neq\Delta(\CP)$, and let $m=m(S)$. The following conditions are equivalent:
\begin{enumerate}
\item $\CP\rk(S)=1$.
\item $w-m\in\langle\mathcal P\rangle$ for every $w\in\Ap(S,m)$ with $0<w<\max(\mathcal P)$.
\end{enumerate}
Moreover, if these conditions hold, then $S$ has maximal embedding dimension, $t(S)=m-1$ and $\msg(S)=\{m\}\cup(\{m\}+\PF(S))$.
\end{theorem}

\begin{proof}
(1)$\Rightarrow$(2). By Proposition \ref{proposcion T rango uno}, $S=\{0\}\cup(\{m\}+\langle\mathcal P\cup\{m\}\rangle)\cup\{\max(\mathcal P)+1,\rightarrow\}$. Let $w\in\Ap(S,m)$ with $0<w<\max(\mathcal P)$. Then $w-m\in\langle\mathcal P\cup\{m\}\rangle$. Note that $w\neq m$, because $m\notin\Ap(S,m)$. If $m$ appeared in some expression of $w-m$ as a combination of $\mathcal P\cup\{m\}$, then $w-m\in\{m\}+\langle\mathcal P\cup\{m\}\rangle\subseteq S$, which is impossible. Hence $w-m\in\langle\mathcal P\rangle$.

(2)$\Rightarrow$(1). We prove by induction on $s$ that every $s\in S$ with $0<s<\max(\mathcal P)$ belongs to $\{m\}+\langle\mathcal P\cup\{m\}\rangle$. If $s\in\Ap(S,m)$, then $s-m\in\langle\mathcal P\rangle$ by hypothesis. Otherwise $s-m\in S$; if $s=m$ there is nothing to prove, and if $s-m\neq0$, then by the induction hypothesis $s-m\in\{m\}+\langle\mathcal P\cup\{m\}\rangle$, and so $s\in\{m\}+\langle\mathcal P\cup\{m\}\rangle$. Therefore $S\subseteq T=\{0\}\cup(\{m\}+\langle\mathcal P\cup\{m\}\rangle)\cup\{\max(\mathcal P)+1,\rightarrow\}$. By Lemma \ref{lema M contenido en S}, $T\subseteq S$, so $S=T$. Since $S\neq\Delta(\CP)$, we have $m<\max(\mathcal P)$, and $(\{m\}+\langle\mathcal P\cup\{m\}\rangle)\cap\mathcal P\subseteq S\cap\mathcal P=\emptyset$. By Proposition \ref{proposcion T rango uno}, $\CP\rk(S)=1$.

Now assume that (1) and (2) hold, and let $w\in\Ap(S,m)\setminus\{0\}$. We show that $w-m\in\PF(S)$. Note that $w\neq\max(\mathcal P)$ because $\max(\mathcal P)\notin S$.
\begin{itemize}
\item If $w>\max(\mathcal P)$, then $w-m\notin S$ and, for every $s\in S\setminus\{0\}$, $w-m+s\geq w>\max(\mathcal P)=F(S)$, so $w-m+s\in S$. Hence $w-m\in\PF(S)$.
\item If $w<\max(\mathcal P)$, then $w-m\in\langle\mathcal P\rangle\subseteq\langle\PF(S)\rangle$ by (2), $w-m\neq0$ and $w-m\notin S$. By Lemma \ref{lema sumas de PF}, $w-m\in\PF(S)$.
\end{itemize}
Hence $\{w-m\mid w\in\Ap(S,m)\setminus\{0\}\}\subseteq\PF(S)$. This set has $m-1$ elements and $t(S)\leq m-1$, so $t(S)=m-1$, and Lemma \ref{lema MED tipo} shows that $S$ has maximal embedding dimension and $\Ap(S,m)\setminus\{0\}=\{m\}+\PF(S)$. Finally, for a MED-semigroup, $\msg(S)=\{m\}\cup(\Ap(S,m)\setminus\{0\})$.
\end{proof}

For a MED-semigroup the description $\Ap(S,m)=\{0\}\cup(\{m\}+\PF(S))$ is classical, and from it the following result is immediate. The content of Theorem \ref{teorema MED} is the converse: rank one in $\CP$ forces maximal embedding dimension.

\begin{corollary}\label{corolario MED}
Let $S$ be a numerical semigroup with $m(S)<F(S)$. The following conditions are equivalent:
\begin{enumerate}
\item $S$ has maximal embedding dimension.
\item $\mathcal C(\PF(S))\rk(S)=1$.
\item $S=\{0\}\cup\big(\{m(S)\}+\langle\PF(S)\cup\{m(S)\}\rangle\big)\cup\{F(S)+1,\rightarrow\}$.
\end{enumerate}
\end{corollary}

\begin{proof}
Note that $S\in\mathcal C(\PF(S))$ and that $S\neq\Delta(\mathcal C(\PF(S)))=\{0,F(S)+1,\rightarrow\}$ because $m(S)<F(S)$.

(1)$\Rightarrow$(2). By Lemma \ref{lema MED tipo}, $w-m(S)\in\PF(S)$ for every $w\in\Ap(S,m(S))\setminus\{0\}$, so condition (2) of Theorem \ref{teorema MED} holds with $\mathcal P=\PF(S)$.

(2)$\Rightarrow$(1). It is Theorem \ref{teorema MED}.

(2)$\Leftrightarrow$(3). It is Proposition \ref{proposcion T rango uno} (and Lemma \ref{lema caracterizancion rango 1}) with $\mathcal P=\PF(S)$.
\end{proof}

If $m(S)>F(S)$, then $S=\{0,F(S)+1,\rightarrow\}$, which has maximal embedding dimension and is the minimum of $\mathcal C(\PF(S))$, so it has rank zero. Therefore, apart from these trivial semigroups, the MED-semigroups are exactly the semigroups that are generated by their multiplicity inside the covariety determined by their set of pseudo-Frobenius numbers.

Recall that $n(S)$ denotes the number of elements of $S$ smaller than $F(S)$. Wilf \cite{wilf} asked whether every numerical semigroup $S$ satisfies
\[
F(S)+1\leq e(S)\,n(S),
\]
and this question is known as Wilf's conjecture (see the survey \cite{delgado-wilf}). It is well known that MED-semigroups satisfy it; this is a particular case of the result of Sammartano \cite{sammartano}, who proved the inequality for all the numerical semigroups with $2e(S)\geq m(S)$. For a MED-semigroup it also follows from \eqref{eq FGH}, since $t(S)+1=m(S)=e(S)$ by Lemma \ref{lema MED tipo}. Therefore, as a consequence of Theorem \ref{teorema MED}, we obtain the following result.

\begin{corollary}\label{corolario Wilf}
Let $S\in\CP$ with $\CP\rk(S)=1$. Then $S$ satisfies Wilf's conjecture.
\end{corollary}

\begin{example}\label{ejemplo rango uno 57}
In $\mathcal C(\{5,7\})$, the semigroups of rank one are $\langle 4,9,10,11\rangle$, $\langle 6,8,9,10,11,13\rangle$ and $\langle 3,8,10\rangle$ (see Example \ref{ejemplo generadores 57}). Their types are $3$, $5$ and $2$, that is, $m-1$ in each case. The semigroup $S=\langle 4,6,9,11\rangle$ also has maximal embedding dimension, but its $\mathcal C(\{5,7\})$-rank is $2$: here $\Ap(S,4)=\{0,6,9,11\}$ and $6-4=2\notin\langle 5,7\rangle$, so condition (2) of Theorem \ref{teorema MED} fails. Since $\PF(S)=\{2,5,7\}$, Corollary \ref{corolario MED} tells us that $S=\mathcal C(\{2,5,7\})[\{4\}]$; indeed,
\[
\{0\}\cup(\{4\}+\langle 2,4,5,7\rangle)\cup\{8,\rightarrow\}=\{0,4,6,8,\rightarrow\}=S.
\]
\end{example}

\begin{example}\label{ejemplo tabla rango}
Table \ref{tabla rango} gives, for each $S\in\mathcal C(\{9,11\})$, numbered as in Table \ref{tabla 9 11}, its multiplicity, embedding dimension and type, its $\mathcal C(\{9,11\})$-minimal system of generators (by Example \ref{ejemplo msg}, the minimal generators of $S$ smaller than $11$) and its $\mathcal C(\{9,11\})$-rank. The last two columns say whether $S$ has maximal embedding dimension and whether $S$ is a maximal element of $\mathcal C(\{9,11\})$.
 We observe the following.
\begin{enumerate}
\item The semigroups of rank one are $S_1,\ldots,S_6$, and all of them have maximal embedding dimension, as Theorem \ref{teorema MED} predicts.
\item The converse fails inside $\mathcal C(\{9,11\})$: $S_9$, $S_{11}$, $S_{12}$, $S_{13}$, $S_{18}$ and $S_{20}$ have maximal embedding dimension and rank $2$ or $3$. The reason is that their sets of pseudo-Frobenius numbers are larger than $\{9,11\}$, while condition (2) of Theorem \ref{teorema MED} only allows differences in $\langle 9,11\rangle$. For instance, $\Ap(S_{12},8)=\{0,10,12,13,14,15,17,19\}$ and $10-8=2\notin\langle 9,11\rangle$, so $S_{12}$ does not have rank one in $\mathcal C(\{9,11\})$. However, $2\in\PF(S_{12})$, and by Corollary \ref{corolario MED} the semigroup $S_{12}$ has rank one in $\mathcal C(\PF(S_{12}))=\mathcal C(\{2,4,5,6,7,9,11\})$. The same happens with the other five semigroups.
\item Maximal embedding dimension and maximality in $\mathcal C(\{9,11\})$ are independent properties. Among the maximal elements $S_{20}$, $S_{21}$ and $S_{22}$, only $S_{20}$ has maximal embedding dimension, whereas $S_1,\ldots,S_6$ have maximal embedding dimension and are not maximal.
\end{enumerate}
\end{example}

{\small
\setlength{\tabcolsep}{6pt}
\begin{longtable}{@{}cccclccc@{}}
\caption{Ranks, maximal embedding dimension and maximal elements in $\mathcal C(\{9,11\})$.}\label{tabla rango}\\
\toprule
 & $m(S)$ & $e(S)$ & $t(S)$ & $\mathcal C(\mathcal P)\msg(S)$ & $\mathcal C(\mathcal P)\rk(S)$ & MED & maximal\\
\midrule
\endfirsthead
\multicolumn{8}{@{}l}{\emph{Table \ref{tabla rango} (continued)}}\\
\toprule
 & $m(S)$ & $e(S)$ & $t(S)$ & $\mathcal C(\mathcal P)\msg(S)$ & $\mathcal C(\mathcal P)\rk(S)$ & MED & maximal\\
\midrule
\endhead
\bottomrule
\endlastfoot
$S_{0}$ & 12 & 12 & 11 & $\emptyset$ & 0 & yes & no\\
$S_{1}$ & 6 & 6 & 5 & $\{6\}$ & 1 & yes & no\\
$S_{2}$ & 7 & 7 & 6 & $\{7\}$ & 1 & yes & no\\
$S_{3}$ & 8 & 8 & 7 & $\{8\}$ & 1 & yes & no\\
$S_{4}$ & 10 & 10 & 9 & $\{10\}$ & 1 & yes & no\\
$S_{5}$ & 4 & 4 & 3 & $\{4\}$ & 1 & yes & no\\
$S_{6}$ & 5 & 5 & 4 & $\{5\}$ & 1 & yes & no\\
$S_{7}$ & 6 & 5 & 4 & $\{6, 7\}$ & 2 & no & no\\
$S_{8}$ & 6 & 5 & 4 & $\{6, 8\}$ & 2 & no & no\\
$S_{9}$ & 6 & 6 & 5 & $\{6, 10\}$ & 2 & yes & no\\
$S_{10}$ & 7 & 6 & 5 & $\{7, 8\}$ & 2 & no & no\\
$S_{11}$ & 7 & 7 & 6 & $\{7, 10\}$ & 2 & yes & no\\
$S_{12}$ & 8 & 8 & 7 & $\{8, 10\}$ & 2 & yes & no\\
$S_{13}$ & 4 & 4 & 3 & $\{4, 10\}$ & 2 & yes & no\\
$S_{14}$ & 5 & 4 & 3 & $\{5, 7\}$ & 2 & no & no\\
$S_{15}$ & 5 & 4 & 3 & $\{5, 8\}$ & 2 & no & no\\
$S_{16}$ & 6 & 4 & 3 & $\{6, 7, 8\}$ & 3 & no & no\\
$S_{17}$ & 6 & 4 & 3 & $\{6, 7, 10\}$ & 3 & no & no\\
$S_{18}$ & 6 & 6 & 5 & $\{6, 8, 10\}$ & 3 & yes & no\\
$S_{19}$ & 7 & 5 & 4 & $\{7, 8, 10\}$ & 3 & no & no\\
$S_{20}$ & 4 & 4 & 3 & $\{4, 6\}$ & 2 & yes & yes\\
$S_{21}$ & 5 & 3 & 2 & $\{5, 7, 8\}$ & 3 & no & yes\\
$S_{22}$ & 6 & 4 & 2 & $\{6, 7, 8, 10\}$ & 4 & no & yes\\
\end{longtable}
}

\subsection{\texorpdfstring{Computing the $\CP$-semigroups of rank one}{Computing the C(P)-semigroups of rank one}}\label{subsec:computing}

We now aim to provide an algorithm that allows us to compute all $\mathcal C(\mathcal P)$-semigroups of $\mathcal C(\mathcal P)$-rank one. Note that, by Proposition \ref{proposcion T rango uno}, this is equivalent to computing all $m\in\N\setminus\{0,1\}$ such that $m<\max(\mathcal P)$ and $(\{m\}+\langle\mathcal P\cup\{m\}\rangle)\cap\mathcal P=\emptyset$.

We generalize the concept of the Apéry set to submonoids of $(\N,+)$ as follows. If $M$ is a submonoid of $(\N,+)$ and $m\in M\setminus\{0\}$, then
\[
\mathrm{Ap}(M,m)=\{x\in M \mid x-m\notin M\}.
\]
Note that if $d=\gcd(M)$, then $\frac{M}{d}=\{\frac{x}{d}\mid x\in M\}$ is a numerical semigroup. We thus easily deduce the following result.

\begin{lemma}\label{lema cociente monoide y apery}
    If $M$ is a submonoid of $(\N,+)$, $d=\gcd(M)$ and $m\in M\setminus\{0\}$, then $\mathrm{Ap}(M,m)$ has cardinality $\frac{m}{d}$. Moreover,
    \[
    \mathrm{Ap}(M,m)=\{dw \mid w\in \mathrm{Ap}(\tfrac{M}{d},\tfrac{m}{d})\}.
    \]
\end{lemma}

\begin{remark}
    In the literature there exist numerous algorithms (see, for example, \cite{numericalsgps}) that, given a numerical semigroup $S$ and $n\in S\setminus\{0\}$, compute $\mathrm{Ap}(S,n)$. Therefore, as a consequence of Lemma \ref{lema cociente monoide y apery}, we obtain algorithms to compute $\mathrm{Ap}(M,m)$ when $M$ is a submonoid of $(\N,+)$ and $m\in M\setminus\{0\}$.
\end{remark}

\begin{lemma}\label{lema Apery vacio}
    Let $m\in\N\setminus\{0,1\}$. Then $(\{m\}+\langle\mathcal P\cup\{m\}\rangle)\cap \mathcal P=\emptyset$ if and only if $\mathcal P\subseteq \mathrm{Ap}(\langle\mathcal P\cup\{m\}\rangle,m)$.
\end{lemma}

\begin{proof}
    (Necessity) It is clear that $\mathcal P\subseteq\langle \mathcal P\cup\{m\}\rangle$. If $p\in\mathcal P$ and $p-m\in\langle\mathcal P\cup\{m\}\rangle$, then $p\in\{m\}+\langle \mathcal P\cup\{m\} \rangle$, which is a contradiction.

    (Sufficiency) If there exists $p\in\mathcal P$ such that $p\in\{m\}+\langle\mathcal P\cup\{m\}\rangle$, then $p-m\in\langle\mathcal P\cup\{m\}\rangle$, which contradicts $p\in\Ap(\langle\mathcal P\cup\{m\}\rangle,m)$.
\end{proof}

Using Proposition \ref{proposcion T rango uno} and Lemma \ref{lema Apery vacio}, we obtain that the set of $\mathcal C(\mathcal P)$-semigroups of $\mathcal C(\mathcal P)$-rank one is
\[
\{\mathcal C(\mathcal P)[\{m\}] \mid 1<m<\max(\mathcal P),\mathcal P\subseteq\mathrm{Ap}(\langle\mathcal P\cup\{m\}\rangle,m)\}.
\]

Note that since $\mathcal P\subseteq \mathrm{Ap}(\langle\mathcal P\cup\{m\}\rangle,m)$ and $0\notin\mathcal P$, it follows that $\mathrm{Card}(\mathcal P)+1\leq m$.

The condition $\mathcal P\subseteq\Ap(\langle\mathcal P\cup\{m\}\rangle,m)$ can often be decided by looking at residues modulo $m$.

\begin{lemma}\label{lema congruencias}
Let $m\in\N\setminus\{0,1\}$.
\begin{enumerate}
\item If $\mathcal P\subseteq\Ap(\langle\mathcal P\cup\{m\}\rangle,m)$, then no element of $\mathcal P$ is a multiple of $m$, and the elements of $\mathcal P$ are pairwise incongruent modulo $m$.
\item If $m>\max(\mathcal P)-\min(\mathcal P)$, then $\mathcal P\subseteq\Ap(\langle\mathcal P\cup\{m\}\rangle,m)$ if and only if no element of $\mathcal P$ is a multiple of $m$.
\end{enumerate}
\end{lemma}

\begin{proof}
(1) Let $M=\langle\mathcal P\cup\{m\}\rangle$. If $p=km$ with $p\in\mathcal P$, then $k\geq1$ and $p-m=(k-1)m\in M$. If $p<q$ are elements of $\mathcal P$ with $q=p+km$, then $k\geq 1$ and $q-m=p+(k-1)m\in M$. In both cases we contradict $\mathcal P\subseteq\Ap(M,m)$.

(2) By (1), it is enough to prove the sufficiency. Let $p\in\mathcal P$ and assume that $p-m\in M$. Since $p-m\leq\max(\mathcal P)-m<\min(\mathcal P)$, no element of $\mathcal P$ can appear in an expression of $p-m$, so $p-m=km$ for some $k\in\N$, and $m$ divides $p$, a contradiction.
\end{proof}

\begin{algo}{Computing the $\mathcal C(\mathcal P)$-semigroups with $\mathcal C(\mathcal P)$-rank one}{algoritmo segundo}
\KwIn{A finite nonempty subset $\mathcal P$ of $\mathbb N\setminus\{0\}$.}
\KwOut{The set $R(\mathcal P)=\{m\in\N\mid \mathcal C(\mathcal P)[\{m\}] \text{ is a }\mathcal C(\mathcal P)\text{-semigroup of }\mathcal C(\mathcal P)\text{-rank one}\}$, together with these semigroups.}
$F\gets\max(\mathcal P)$, \quad $R\gets\emptyset$\;
\For{$m=\mathrm{Card}(\mathcal P)+1$ \KwTo $F-1$}{
  \If{$m\mid p$ for some $p\in\mathcal P$ \KwSty{or} $p\equiv q \pmod{m}$ for some distinct $p,q\in\mathcal P$}{
    \KwSty{continue}\Comment*[r]{Lemma \ref{lema congruencias}(1)}
  }
  \uIf{$m>F-\min(\mathcal P)$}{
    $R\gets R\cup\{m\}$\Comment*[r]{Lemma \ref{lema congruencias}(2)}
  }
  \ElseIf{$\mathcal P\subseteq\Ap(\langle\mathcal P\cup\{m\}\rangle,m)$}{
    $R\gets R\cup\{m\}$\Comment*[r]{Lemma \ref{lema Apery vacio}}
  }
}
\KwRet{$R$ and, for each $m\in R$, $\mathcal C(\mathcal P)[\{m\}]=\{0\}\cup(\{m\}+\langle\mathcal P\cup\{m\}\rangle)\cup\{F+1,\rightarrow\}$}\;
\end{algo}

By Theorem \ref{teorema MED}, every semigroup returned by Algorithm \ref{algoritmo segundo} has maximal embedding dimension and type $m-1$. In particular, Algorithm \ref{algoritmo segundo} produces, without building the tree, semigroups of $\CP$ of every type $m-1$ with $m\in R(\mathcal P)$.

\begin{example}
    If we take $\mathcal P=\{5,7\}$, then Algorithm \ref{algoritmo segundo} yields that
    \[
    \{m\in\N\mid \mathcal C(\{5,7\})[\{m\}] \text{ is a } \mathcal C(\{5,7\})\text{-semigroup of } \mathcal C(\{5,7\})\text{-rank one}\}=\{3,4,6\}.
    \]
    Here $\max(\mathcal P)-\min(\mathcal P)=2$, so for every candidate $m\in\{3,4,5,6\}$ the decision is taken by Lemma \ref{lema congruencias}: $m=5$ is rejected because $5\in\mathcal P$, and $3$, $4$, $6$ are accepted. The corresponding semigroups are $\langle 3,8,10\rangle$, $\langle 4,9,10,11\rangle$ and $\langle 6,8,9,10,11,13\rangle$, of types $2$, $3$ and $5$.
\end{example}

\begin{example}
Let $\mathcal P=\{9,11\}$. The candidates are $m\in\{3,\ldots,10\}$, and $\max(\mathcal P)-\min(\mathcal P)=2<m$ for all of them, so again Lemma \ref{lema congruencias}(2) decides: $m=3$ and $m=9$ are rejected because they divide $9$, and $R(\{9,11\})=\{4,5,6,7,8,10\}$. The semigroups obtained are
\begin{gather*}
\langle 4,13,14,15\rangle,\ \langle 5,12,13,14,16\rangle,\ \langle 6,13,14,15,16,17\rangle,\\ \langle 7,12,13,15,16,17,18\rangle,\ \langle 8,12,13,14,15,17,18,19\rangle,\ \langle 10,12,13,\ldots,19,21\rangle,
\end{gather*}
that is, $S_5,S_6,S_1,S_2,S_3,S_4$ of Table \ref{tabla 9 11}, with types $3,4,5,6,7,9$. Note that $\mathcal C(\{9,11\})$ has no semigroup of rank one and type $8$, because $9\notin R(\{9,11\})$.
\end{example}

We conclude this section with a way of producing, from the Ap\'ery sets of a numerical semigroup, sets $\mathcal P$ together with $\CP$-semigroups of $\CP$-rank one.

\begin{corollary}\label{corolario algoritmo}
    Let $S$ be a numerical semigroup, $S\ne\N$, and $m\in S\setminus\{0\}$. If $\mathcal P\subseteq\mathrm{Ap}(S,m)\setminus\{0\}$ and $m<\max(\mathcal P)$, then
    \[
    T=\{0\}\cup(\{m\}+\langle\mathcal P\cup\{m\}\rangle)\cup\{\max(\mathcal P)+1,\rightarrow\}
    \]
    is a $\mathcal C(\mathcal P)$-semigroup of $\mathcal C(\mathcal P)$-rank one. Moreover, every $\mathcal C(\mathcal P)$-semigroup of $\mathcal C(\mathcal P)$-rank one can be obtained in this way, for a suitable choice of $S$ and $m$.
\end{corollary}

\begin{proof}
    To prove the first part, using Proposition \ref{proposcion T rango uno} and Lemma \ref{lema Apery vacio}, it suffices to show that $\mathcal P\subseteq\mathrm{Ap}(\langle\mathcal P\cup\{m\}\rangle,m)$. This holds since if $x\in\mathcal P$ and $x-m\in\langle\mathcal P\cup\{m\}\rangle$, then $x-m\in S$, which contradicts the fact that $\mathcal P\subseteq\mathrm{Ap}(S,m)$.

    If $T$ is a $\mathcal C(\mathcal P)$-semigroup of $\mathcal C(\mathcal P)$-rank one, then applying again Proposition \ref{proposcion T rango uno} and Lemma \ref{lema Apery vacio}, we have that
    \[
    T=\{0\}\cup(\{m\}+\langle\mathcal P\cup\{m\}\rangle)\cup\{\max(\mathcal P)+1,\rightarrow\}
    \]
    for some $m\in\N$ such that $1<m<\max(\mathcal P)$ and $\mathcal P\subseteq\mathrm{Ap}(\langle\mathcal P\cup\{m\}\rangle,m)$. The proof is completed by taking $S=\langle\mathcal P\cup\{m\}\rangle\cup\{\max(\mathcal P)+1,\rightarrow\}$; this is a numerical semigroup, $m\in S$ and $\mathcal P\subseteq\Ap(S,m)$, because $x-m<\max(\mathcal P)+1$ for every $x\in\mathcal P$.
\end{proof}

\begin{example}
    Let $S=\langle5,7,9\rangle=\{0,5,7,9,10,12,14,\rightarrow\}$, $5\in S\setminus\{0\}$ and $\mathrm{Ap}(S,5)=\{0,7,9,16,18\}$.
    \begin{itemize}
        \item Since $\{16,18\}\subseteq\mathrm{Ap}(S,5)$ and $5<18$, Corollary \ref{corolario algoritmo} with $\mathcal P=\{16,18\}$ and $m=5$ shows that
        \[
        T=\{0\}\cup(\{5\}+\langle5,16,18\rangle)\cup\{19,\rightarrow\}
        \]
        is a $\mathcal C(\{16,18\})$-semigroup of $\mathcal C(\{16,18\})$-rank one. Moreover, $T=\mathcal C(\{16,18\})[\{5\}]$. Explicitly, $T=\langle 5,19,21,22,23\rangle$; it has maximal embedding dimension and $\PF(T)=\{14,16,17,18\}$, so $t(T)=4=m(T)-1$, in agreement with Theorem \ref{teorema MED}. Note that $\PF(T)$ is strictly larger than $\{16,18\}$.
        \item Similarly, since $\{7,18\}\subseteq\mathrm{Ap}(S,5)$, Corollary \ref{corolario algoritmo} with $\mathcal P=\{7,18\}$ shows that
        \[
        T=\{0\}\cup(\{5\}+\langle5,7,18\rangle)\cup\{19,\rightarrow\}
        \]
        is a $\mathcal C(\{7,18\})$-semigroup of $\mathcal C(\{7,18\})$-rank one. Moreover, $T=\mathcal C(\{7,18\})[\{5\}]$. Explicitly, $T=\langle 5,12,19,21,23\rangle$ and $\PF(T)=\{7,14,16,18\}$.
    \end{itemize}
\end{example}

\noindent{\bf Fundings.} First named author was supported by the Spanish grant PID2023-149508NB-I00, funded by MICIU/AEI /10.13039/501100011033 and by FEDER, UE.

\noindent {\sc William Giovanni Hern\'andez-Yanes}\\
Departamento de Matem\'aticas, Estad\'{\i}stica e Investigaci\'on Operativa\\
	Instituto Universitario de Matem\'aticas y Aplicaciones (IMAULL)\\ Universidad de La Laguna. Apartado de Correos 456.
	38200 La Laguna, Tenerife, Espa\~na \\ORCID ID: 0009-0002-5140-8295
	
	 \noindent {whernany@ull.es}
	\vspace{0.3cm}

\noindent {\sc José Carlos Rosales}\\
Department of Algebra, Faculty of Sciences, University of Granada\\ E-18071 Granada, Spain \\ORCID ID: 0000-0003-3353-4335
	
	 \noindent {jrosales@ugr.es}
	

\begin{thebibliography}{99}

\bibitem{apery} Ap\'ery, R. \textit{Sur les branches superlinéaires des courbes algébriques}. C. R. Acad. Sci. Paris 222 (1946), 1198--1200.

\bibitem{barucci} Barucci, V.; Dobbs, D.E.; Fontana, M. \textit{Maximality Properties in Numerical Semigroups and Applications to One-Dimensional Analytically Irreducible Local Domains}. Memoirs Amer. Math. Soc. 125 (598), 1997.

\bibitem{barucci-froberg} Barucci, V.; Fr\"oberg, R. \textit{One-dimensional almost Gorenstein rings}. J. Algebra 188 (1997), no. 2, 418--442.

\bibitem{delgado-wilf} Delgado, M. \textit{Conjecture of Wilf: a survey}. In: Barucci, V.; Chapman, S.; D'Anna, M.; Fr\"oberg, R. (eds.), \textit{Numerical Semigroups}, Springer INdAM Series 40, Springer, Cham, 2020, 39--62.

\bibitem{numericalsgps} Delgado, M.; Garc\'{\i}a-S\'anchez, P.A. \textit{numericalsgps, a GAP package for numerical semigroups}. ACM Commun. Comput. Algebra 50 (2016), no. 1, 12--24.

\bibitem{numericalsgps-pkg} Delgado, M.; Garc\'{\i}a-S\'anchez, P.A.; Morais, J. \textit{NumericalSgps, a package for numerical semigroups}, Version 1.5.0 (2026). \texttt{GAP} package, \url{https://gap-packages.github.io/numericalsgps}.

\bibitem{delgado} Delgado, M.; Garc\'{\i}a-S\'anchez, P.A.; Robles-P\'erez, A.M. \textit{Numerical semigroups with a given set of pseudo-Frobenius numbers}. LMS J. Comput. Math. 19 (2016), no. 1, 186--205.

\bibitem{froberg} Fr\"oberg, R.; Gottlieb, C.; H\"aggkvist, R. \textit{On Numerical Semigroups}. Semigroup Forum 35 (1987), 63--83.

\bibitem{covariedades} Moreno-Fr\'{\i}as, M.A.; Rosales, J.C. \textit{The Covariety of Numerical Semigroups With Fixed Frobenius Number}. J. Algebraic Combin. 60 (2024), no. 2, 555--568.

\bibitem{nari} Nari, H. \textit{Symmetries on almost symmetric numerical semigroups}. Semigroup Forum 86 (2013), 140--154.

\bibitem{JPAA} Rosales, J.C.; Branco, M.B. \textit{Numerical Semigroups that can be expressed as an Intersection of Symmetric Numerical Semigroups}, J. Pure Appl. Algebra 171 (2002), 303--314.

\bibitem{libro} Rosales, J.C.; Garc\'{\i}a-S\'anchez, P.A. \textit{Numerical Semigroups}. Developments in Mathematics, Vol. 20, Springer, New York, 2009.

\bibitem{med} Rosales, J.C.; Garc\'{\i}a-S\'anchez, P.A.; Garc\'{\i}a-Garc\'{\i}a, J.I.; Branco, M.B. \textit{Numerical semigroups with maximal embedding dimension}. In: Badawi, A. (ed.), \textit{Focus on Commutative Rings Research}, Nova Science Publishers, New York, 2006, 47--53.

\bibitem{sammartano} Sammartano, A. \textit{Numerical semigroups with large embedding dimension satisfy Wilf's conjecture}. Semigroup Forum 85 (2012), no. 3, 439--447.

\bibitem{wilf} Wilf, H.S. \textit{A circle-of-lights algorithm for the ``money-changing problem''}. Amer. Math. Monthly 85 (1978), 562--565.

\end{thebibliography}
\end{document}